\documentclass[11pt]{article}

\usepackage[letterpaper,margin=1in]{geometry}
\usepackage{amssymb,amsmath,amsfonts,eurosym,geometry,ulem,graphicx, color,sectsty,setspace, comment,footmisc,pdflscape,array}
\usepackage{subcaption}
\definecolor{bur}{RGB}{158,5,8}

\newcolumntype{L}[1]{>{\raggedright\let\newline\\arraybackslash\hspace{0pt}}m{#1}}
\newcolumntype{C}[1]{>{\centering\let\newline\\arraybackslash\hspace{0pt}}m{#1}}
\newcolumntype{R}[1]{>{\raggedleft\let\newline\\arraybackslash\hspace{0pt}}m{#1}}
\usepackage{mathtools}

\usepackage{xcolor}
\usepackage{ulem}
\usepackage{graphicx}
\usepackage{enumitem}
\usepackage{algorithm}
\usepackage{algpseudocode}
\usepackage{bbm}

\usepackage[para,online,flushleft]{threeparttable}
\usepackage{booktabs,multirow,makecell,tabularx}
\definecolor{darkblue}{rgb}{0.0, 0.0, 0.55}
\definecolor{lm}{rgb}{0.95, 0.4, 0.8}
\usepackage[colorlinks=true,linkcolor=blue!80!black,citecolor=blue!80!black,urlcolor=red!80!black]{hyperref}
\usepackage{longtable}

\usepackage{bbm}

\usepackage{booktabs}

\def\min{\mathop{\rm min}}
\def\max{\mathop{\rm max}}

\usepackage{xcolor}

\newif \iffinal
\finalfalse
\iffinal
\newcommand{\EDmodified}[1]{{#1}}
\newcommand{\EDcomments}[1]{{}}
\else
\newcommand{\EDmodified}[1]{{\color{blue} #1}}
\newcommand{\EDcomments}[1]{{\EDmodified{Erick commented: #1}}}

\fi
\usepackage{tikz}
\usepackage{pgfplots}
\pgfplotsset{compat=1.10}

\usepackage{algorithm}
\usepackage{algpseudocode}
\usepackage{graphicx}
\usepackage{caption}
\usepackage{subcaption}
\usepackage{float}
\definecolor{lm}{rgb}{0.95, 0.4, 0.8}
\newcommand{\USmodified}[1]{{#1}}

\usepackage{graphicx}
\usepackage{float}
\usepackage{subcaption}
\usepackage{bbm}

\newtheorem{theorem}{Theorem}
\newtheorem{definition}{Definition}

\newtheorem{proposition}{Proposition} %
\newenvironment{proof}[1][Proof]{\noindent\textbf{#1.} }{\ \rule{0.5em}{0.5em}}

\usepackage[margin = 1.5cm, font = small, labelfont = bf, labelsep = period]{caption}
\usepackage{natbib}
\begin{document}
\title{Subsidizing a New Technology: An Impulse Stackelberg Game Approach%
\thanks{%
Research supported by NSERC, Canada, grant {RGPIN-2024-05067} (Utsav
Sadana) and grant RGPIN-2021-02462 (Georges Zaccour).}}
\author{Utsav Sadana \\
CIRRELT, GERAD and Department of Computer Science\\ and Operations Research, Universit\'{e} de Montr\'{e}al, Canada\thanks{%
utsav.sadana@umontreal.ca} \and Georges Zaccour \\
GERAD and HEC Montr\'{e}al, Canada\thanks{%
georges.zaccour@gerad.ca}}
\maketitle

\begin{abstract}
Governments are motivated to subsidize profit-driven firms that manufacture
zero-emission vehicles to ensure they become price-competitive. This paper
introduces a dynamic Stackelberg game to determine the government's optimal
subsidy strategy for zero-emission vehicles, taking into account the pricing
decisions of a profit-maximizing firm. While firms have the flexibility to
change prices continuously, subsidies are adjusted at specific time
intervals. This is captured in our game formulation by using impulse
controls for discrete-time interventions. We provide a verification theorem
to characterize the Feedback Stackelberg equilibrium and illustrate our
results with numerical experiments.

\end{abstract}

\section{Introduction}

To reduce the greenhouse gas emissions, governments around the world are
offering subsidies to encourage consumers to buy electric vehicles instead
of gasoline cars. In US, the subsidy takes the form of tax credits that top
out at \$7,500 in 2024.\footnote{{{\href{https://www.edmunds.com/fuel-economy/the-ins-and-outs-of-electric-vehicle-tax-credits.html}%
{%
https://www.edmunds.com/fuel-economy/the-ins-and-outs-of-electric-vehicle-tax-credits.html%
}}}} European countries are also offering incentives to consumers to buy
electric vehicles. To illustrate, France gives a subsidy of \euro 5,000,
Italy \euro 3,000, The Netherlands \euro 2,950, and Spain up to \euro 7,000,
plus some other benefits to consumers buying a BEV (Battery Electric
Vehicle).\footnote{{{\href{https://www.fleeteurope.com/en/new-energies/europe/features/ev-incentives-2024-europes-major-fleet-markets?a=FJA05\&t\%5B0\%5D=Taxation\&t\%5B1\%5D=EVs\&curl=1}%
{%
https://www.fleeteurope.com/en/new-energies/europe/features/ev-incentives-2024-europes-major-fleet-markets?a=FJA05\&t\%5B0\%5D=Taxation\&t\%5B1\%5D=EVs\&curl=1%
}}}} Subsidy programs are (normally) designed with a target in mind, and
typically have an end date. For instance, the Canadian Zero-Emission
Vehicles (ZEV) program target is 100\% new light-weight ZEVs sales by 2035,
and it will run until March 31, 2025, or until available funding is
exhausted.\footnote{{{\href{https://tc.canada.ca/en/road-transportation/innovative-technologies/zero-emission-vehicles/light-duty-zero-emission-vehicles}%
{%
https://tc.canada.ca/en/road-transportation/innovative-technologies/zero-emission-vehicles/light-duty-zero-emission-vehicles%
}}}} Another example is the target set by President Obama in 2011 of
\textquotedblleft one million electric vehicles on the road by
2015.\textquotedblright

\USmodified{Subsidy programs to incentivize consumers to buy a new car
have been around for some time. Following the global recession in 2008, many
countries implemented so-called car scrappage incentive programs in which
consumers get a subsidy when trading in their old car with a new one. In
U.S., an individual could get up to \$4,500 when exchanging their car for a
more fuel-efficient one. The dual objective of a scrappage program is to
stimulate the sales of new cars and reduce emissions as new cars pollute
less than old ones per traveled mile. These programs involved significant
amounts, with the German program having the highest budget of 2,5 billion
euros \citep{Ewing2009}. 

How to measure the success of a car subsidy program? Does the
fact that the \$1 billion U.S. budget was used up well ahead of the expected
date \citep{Huang2014} constitute a measure of success? It surely does,
but other metrics should also be considered. Several contributions aimed at
assessing the efficiency of these programs were published; see, \cite{VanWee2011} for a review. A general result seems to be that these programs
succeeded in boosting sales of new cars, but did not necessarily lead to
significantly less pollution. One reason is that individuals drive longer
distances with their new car than what they used to do with their old one.
In terms of targeting the right consumers, the literature highlights the
heavy presence of free riders, that is, consumers who would have bought
anyway a new car. \cite{Hoekstra2014} estimated the proportion of free
riders to 60\% in the case of US. In the specific case of electric cars, a
recent study showed that electric vehicle tax credits boosted sales of U.S.
automakers, reduced pollution, but mainly benefited consumers who would have
purchased EVs without subsidies.\footnote{\href{https://news.stanford.edu/stories/2024/10/electric-vehicle-subsidies-help-the-climate-u-s-automakers}{https://news.stanford.edu/stories/2024/10/electric-vehicle-subsidies-help-the-climate-u-s-automakers}}

Another important issue, which is directly related to our
work, is the behavior of the manufacturers and car dealers when a subsidy to
consumers is available. \cite{Jimenez2016} found that the average
vehicle price in Spain increased by \EUR600 after the implementation of a
scrappage program. In the case of Germany where consumers\ were offered \EUR%
2,500 when scrapping their old car, \cite{Kaul2016} obtained the counter
intuitive result that subsidized buyers paid a little more than comparable
non-subsidized buyers. (To be subsidized, the age of the old car must be
larger than a pre-determined level. Therefore, not all consumers can benefit
from a trade-in program.) Considering a two-player game with a manufacturer
and a government as players, \cite{Janssens2014} and \cite{zaman2024subsidies} reached the conclusion that the benefits of a subsidy program
are mainly cashed by the manufacturer, not the consumers. This leads to the
obvious conclusion that the strategic behavior of manufacturers should not
be ignored when designing a subsidy program. }

Subsidizing a new durable product aims at achieving a series of objectives,
among them reducing the unit production cost and increasing consumers'
confidence in the product. Indeed, it is empirically documented that the
unit production cost decreases with experience, which is measured by
cumulative production \citep{Levitt2013}. By boosting demand, subsidies
accelerate the drop in the marginal production cost, which in turn should
lead to lower price and higher adoption rate. Further, early adopters of a
new product influence non-adopters purchasing behavior through product
reviews and word-of-mouth communications. By increasing the number of early
adopters, subsidies amplify the social impact that adopters have on non
adopters. Seeing other consumers buying an EV increases awareness of the
product and eventually decreases the perceived risk of adopting this new
type of vehicle. \USmodified{A large body of research explores how subsidies,
demographics, socioeconomic factors, and psychological dimensions shape
consumers' willingness to pay for electric vehicles, to inform effective
government intervention policies 
\citep{helveston2015will,
liao2017consumer, fevang2021goes, muehlegger2022subsidizing}}.

Based on the above discussion, our objective is to answer the following
research questions:

\begin{enumerate}
\item What are the equilibrium price and subsidy and how do they evolve over
time?\label{q1}
\item Does the seller take advantage of the subsidy program to raise its
price? \label{q2}
\item What are the cost and benefit of the subsidy program? \label{q3}
\item What is the effect of varying the parameter values on the results? %
\label{q4}
\end{enumerate}

To answer these questions, we develop a game model with two players, a firm
selling EVs and a government subsidizing consumers when purchasing one. The
objective of the government is to reach a cumulative adoption target with a
minimum budget, whereas the firm maximizes its discounted profit over its
planning horizon. By retaining a game model, we account for the strategic
interactions between the pricing policy of the firm and the subsidy policy
of the government. Further, as the learning-by-doing \USmodified{%
\citep{wright1936factors}} in production and the diffusion effect that
adopters exert on non-adopters are inherently dynamic, so is our model. 
\USmodified{\cite{wright1936factors} shows that airplane manufacturing costs
decrease at a constant rate with each doubling of production. Furthermore,
there is a large body of research on estimating the learning-by-doing
effects on cost of production 
\citep{argote1990learning, ziegler2021re,
sweeney2024thomas}. \cite{barwick2025drive} quantifies how learning-by-doing
reduces EV battery costs, estimating a 7.5 learning rate (i.e., doubling
production cuts costs by 7.5). Furthermore they show that learning-by-doing
massively magnifies EV subsidies' effects, boosting global EV sales by 170\%
(versus 29.9\% without learning-by-doing), showing that ignoring learning-by-doing underestimates policy effectiveness.}

The paper is organized as follows: {Section \ref{section:lit_Review}
reviews the related literature on subsidies and differential games}. Section %
\ref{sec:model} introduces a two-player game between a profit-maximizing
firm that solves linear-quadratic regulator-type problem and a government
that aims to reach a desired adoption rate of ZEVs with minimal budget. In
Section \ref{sec:fse}, we derive the sufficient conditions for
characterizing the {Feedback Stackelberg equilibrium (FSE)}. Numerical
results are {presented} in Section \ref{sec:numerical} and conclusions are
given in Section \ref{sec:conclusion}.

\section{Literature Review}

\label{section:lit_Review}

Our paper belongs to the literature on new durable product diffusion
initiated in the seminal paper by \cite{Bass1969new}.\footnote{%
In 2004, \cite{Bass1969new} was voted one of the ten most influential papers
published in Management Science during the last fifty years.} The early
contributions were forecast oriented, that is, they estimated the parameter
values of the diffusion dynamics equation to predict the adoption rate.%
\footnote{%
This literature typically assumes that consumer buys at most one unit. \cite%
{islam2000modelling} extend this class of forecasting models to account for
replacement purchase of the durable product.} In this stream, the firm is
passive, i.e., it does not make pricing or any other decision. \cite%
{RobinsonLakhani1975} extended the framework to a continuous-time
optimal-control problem where the firm decides on the price at each instant
of time. \USmodified{Kalish (1983) introduced dynamics in both the demand
function (saturation and word-of-mouth effects) and in the cost function
(learning-by-doing effect).} \cite{Eliashberg1986} consider a two-stage
model with a monopoly period followed by a duopoly period, and analyze the
pricing strategies of the incumbent. \cite{DocknerJorgensen1988} introduced
price competition in a dynamic oligopoly using a differential game approach.
Each of these papers were followed by a large number of studies considering
some variations. We shall refrain from reviewing the huge literature on
diffusion models and refer the reader to the surveys and tutorials in \cite{Mahajan1993, Mahajan2000, Jorgensen2004,Peres2010}.
Here, we focus on diffusion models with price subsidy.

\cite{KalishLilien1983} were the first to investigate the effect of price
subsidy on the rate of adoption of an innovation in a new product diffusion
framework. The decision maker is the government that chooses the subsidy
rate to maximize the total number of units sold by the terminal date of the
subsidy program. The industry is assumed to be competitive and does not
behave strategically. \cite{Lilien1984} applies the theory developed in \cite%
{KalishLilien1983} to the National Photovoltaic Program implemented by the
Department of Energy in the United States in the 1970s.

Assuming that the new technology is patented, which prevents entry in the
industry at least in the short run, \cite{Zaccour1996} proposed a
differential game played by a firm and a government. The firm chooses the
price and government {sets} the (varying) subsidy rate over time and an
open-loop Nash equilibrium is determined. As in \cite{KalishLilien1983}, the
objective of the government is to maximize the cumulative sales by the
terminal date of the subsidy, which is assumed to also be the firm's
planning horizon. Under similar assumptions, \cite{Dockner1996} consider the
government to be leader and the firm follower in a Stackelberg game. The
authors characterized and compared open-loop and feedback Stackelberg
pricing and subsidy equilibrium strategies.

\cite{Jorgensen1999} retained the same sequential move structure in \cite%
{Dockner1996} and analyzed open-loop Stackelberg equilibrium in a setup
where the government subsidizes consumers and also purchases some quantity
of the new technology to equip its institutions. Both instruments have the
same objective of accelerating the decrease in the unit production cost
through learning-by-doing. As an open-loop Stackelberg equilibrium is in
general time inconsistent, its implementation requires that the leader will
indeed commit to its announcement. \cite{DeCesare2001} introduced
advertising in a Stackelberg differential game played by a firm and a
government. By doing so, the sales rate is affected by both costless
word-of-mouth communication that emanates from within the social system, and
costly advertising paid for by the firm.

\cite{Janssens2014} criticized the above cited papers on three
grounds. First, they state that there is no empirical support to the
assumption that both players have the same planning horizon. A subsidy
program is short-lived, whereas a firm hopes to remain active in the long
term. Second, the assumption that the unit cost decreases linearly in
cumulative production is also questionable empirically. Finally, maximizing
the number of units sold by a certain date is not the best objective a
government can choose because it could be very costly and does not
necessarily help in bringing down the price after the subsidy program.
Consequently, the authors instead minimize the government's budget needed to
reach a certain target. One drawback in this paper is the use of open-loop
information structure in a Stackelberg game, which leads in general to
time-inconsistent equilibrium strategies.

Should the subsidy be increasing, decreasing, or constant over time? \USmodified{(It is
easy to rule out on economic {grounds} that the subsidy cannot be
non-monotone over time.) This question has been considered in papers
where the game is between the government seeking an optimal subsidy policy
and faresighted consumers, that is, they choose the timing of adoption that
maximizes their intertemporal utility and not only current one; see, e.g.,
Zaman and Zaccour (2021) and Langer and Lemoine (2022). Zaman and Zaccour
(2023) extended their two-stage model by adding a manufacturer that acts
strategically and assess the cost of the subsidy and its impact on consumers
and the manufacturer's profit. In particular, they obtained that the price
is increasing over time and is higher than in the absence of subsidy, which
is clearly not a desirable outcome of a subsidy program. }The determination
of subsidy path depends on the diffusion effect (word of mouth and possibly
saturation effects) and cost dynamics (learning in production). The
assumption retained in the literature is that the government can change
continuously the subsidy over time. Such assumption, which is clearly
motivated by mathematical tractability, is very hard to justify (and
implement) in practice. Government agencies do not have the agility to
continuously change their decisions and if they do have, it is politically
and practically difficult to implement/justify. Indeed, think of a subsidy
that takes the form of a tax credit as in the US, and the government is
changing continuously its level. Further, it is intuitive to assume that
modifying the subsidy level entails a fixed cost that should be normally
considered in the design of the program. \USmodified{A separate literature
explores how subsidies affect investment decisions under market uncertainty,
but typically treats subsidies as exogenously determined rather than
strategically adjusted at discrete points in time 
\citep{bigerna2019green,
barbosa2022optimal,  sendstad2022impact, hagspiel2025green}. In this paper,
by contrast, we focus on dynamic interactions between subsidy and diffusion
models to speed the production of ZEVs.}

In this paper, we depart from the literature and suppose that the government
makes subsidy adjustments at specific dates to reach a desired adoption
target with minimum public spending, while considering the pricing decisions
of the profit-maximizing firm. We retain the assumption that the firm can
change continuously its price and that the game is played \`{a} la
Stackelberg, with the government acting as leader and the firm as follower.
We adopt a feedback-information structure and determine feedback-Stackelberg
equilibrium, which is subgame perfect; see, e.g., \cite{basar1998dynamic}, 
\cite{HKZ2012}, \cite{basar2018handbook} for a discussion of the different
information structures in differential games and resulting equilibria. For
applications of Stackelberg equilibrium in the operations management and
supply chain literature, see the surveys in \cite{He2007} and \cite%
{LiSethi2017}.

The theory of dynamic games has been developed assuming that all players
intervene at all decision moments in the game, that is, continuously in a
differential game and at discrete instants of time in a multistage game. It
is only very recently that some advancements have been made on nonzero-sum
impulse games to study discrete-time interventions in continuous-time
systems 
\citep{Aid2020, basei2021nonzero, Sadana_Reddy_Basar_Zaccour_2021,
Sadana_Reddy_Zaccour_2021,  Sadana_Reddy_Zaccour_2023}. However, these
papers consider Nash equilibrium where players decide on their strategies
simultaneously without knowing the strategy of each other. In our subsidy
model, the dominant view is that a Stackelberg equilibrium should be sought
as quite naturally the government has the option of announcing its strategy
before the firm acts. Consequently, we introduce here a new framework, to
which we shall refer as impulse dynamic Stackelberg game (iDSG), which
incorporates subsidies that are adjusted based on the adoption rate of ZEVs
at discrete instants of time. This approach contrasts with all the papers in
this literature that analyzed a continuous-time dynamic Stackelberg game
(DSG) with continuous control for subsidies, and further distinguishes our
work by assuming that subsidies can only take on discrete values. Again, it
is hard to believe that the subsidy is a continuous variable and having a
discrete variable is more realistic. Furthermore, we provide a verification
theorem to characterize the FSE strategies of the government and the firm
and illustrate our results using numerical experiments.

To {wrap} up, we make two important contributions in this paper. First, by
letting the government {intervene} at discrete moments in time, assuming
discrete values for subsidy adjustments, and having a fixed cost attached to
each adjustment, we believe that our modeling of the strategic interactions
involved in a subsidy program is more realistic than what has been done
before in the literature. Second, to the best of our knowledge, it is the
first paper to characterize the equilibrium in an impulse dynamic
Stackelberg game. This is clearly a significant contribution to the theory
of differential games that opens the door to many potential applications.

\section{Model}

\label{sec:model} In this section, we introduce the two-player Stackelberg
game between the government and the firm. The two players use different
kinds of strategies to influence the cumulative sales of the ZEVs. Whereas
the firm can continuously change the price over time, the government chooses
the subsidy levels only at certain discrete decision dates, $\tau _{1},\tau
_{2},\ldots,\tau _{N}$, where $0\leq \tau _{i}\leq T,\,i=\{0,1,\ldots ,N\}$%
. 

Denote by $p(t)$ the price of a ZEV and by $p_{a}$ the given
price of vehicles using old technology, e.g., gasoline motor. For
simplicity, we assume that $p_{a}$ remains constant throughout the planning
horizon. (Letting $p_{a}$ be defined by a function of time would cause no
conceptual difficulty.) 
The discrete set of subsidies that could be offered to the customers is
denoted by $\mathcal{S}=\{0,s_{1},s_{2}, \USmodified{\ldots}, s_{m}\}$, where $s_{i}>0$
and $0$ corresponds to the case with no subsidy. We let the sales rate of
ZEVs be given by 
\begin{equation}
\dot{x}(t)=\alpha _{1}+\alpha _{2}x(t)-\beta (p(t)-p_{a}),\text{ \ \ }%
x(0)=x_{0},  \label{Demand}
\end{equation}%
where $x(t)$ denotes the cumulative sales at time $t$, and $\alpha
_{1},\alpha _{2}$, and $\beta $ are positive parameters. As in \cite%
{Jorgensen1999} and \cite{Jorgensen2004}, our sales function is linear in
the difference in prices of the two technologies, i.e., $\Delta (t)=$\ $%
p(t)-p_{a}$, and is increasing in $p_{a}$ and decreasing in $p(t) $. To have
non-negative demand, we suppose that $\Delta (t)\leq \frac{\alpha
_{1}+\alpha _{2}x(t)}{\beta }$. The market size, which is given by $\alpha
_{1}+\alpha _{2}x(t)$, is not constant, but endogenous and increasing in
cumulative sales. In Bass's seminal paper \citep{Bass1969new}, the term $%
\alpha _{2}x(t)$ is defined as the word-of-mouth effect, i.e., the positive
impact exerted by adopters on not yet adopters of the new product.
Alternatively to this information dissemination (or free advertising)
interpretation, one can assume that the larger $x(t)$, the easier is to find
a public place to recharge the battery, which in turn enlarges the market
potential and demand.

To decrease the price gap\ between the two technologies and thereby boost
the demand, the government gives a subsidy $s(t)$. Consequently, the demand
becomes 
\begin{equation}
\dot{x}(t)=\alpha _{1}+\alpha _{2}x(t)-\beta (p(t)-s(t)-p_{a}).
\label{eq:state_dyn}
\end{equation}%
If the government changes the subsidy level, at a decision date $\tau $, $%
x(t)$ has a kink at $t=\tau $. We assume that the unit production cost is
decreasing in cumulative sales, which captures the idea of learning-by-doing
effect, and is given by 
\begin{equation*}
c(x(t))=b_{1}-b_{2}x(t),
\end{equation*}%
where $b_{1}$ is the initial unit cost and $b_{2}>0$ measures the learning
speed. We will insure that the cost remains always positive. \USmodified{%
Learning-by-doing, which essentially means a reduction in unit production
cost, has a long history in the economics literature. Indeed, it has been
already mentioned one century ago in \cite{Thorndike1927}, with its economic
implications well documented and discussed over several decades; see, e.g., 
\cite{wright1936factors}, \cite{Arrow1962}, \cite{DuttonThomas1984}, and 
\cite{McDonald2001}.  \cite{Yelle1979} gives a historical review of the
learning-by-doing phenomenon. The assumption that the cost function is
linear in $x(t)$ has also been adopted in, e.g., \cite{fudenberg1983capital}, 
\cite{raman1995optimal}, \cite{Jorgensen2000}, \cite{xu2011dynamic} and \cite%
{creti2018defining}. \cite{Kogan2016} argue that learning is linear in
mature industry, meaning that after a first (possibly short) period of rapid
decrease in the unit production cost, the gain flattened as time goes by.}

The objective of the firm is to maximize its discounted stream of profit
over the planning horizon $T$, that is, 
\begin{equation}
\USmodified{J^{f}(p(\cdot ),\eta (\cdot ))}=\max_{p(t)\in \Omega
^{f}}\int_{0}^{T}e^{-\rho t}(p(t)-c(x(t))\USmodified{)}\dot{x}(t)dt,
\end{equation}%
where $\rho $ is the discount factor, $\eta (\cdot )$ is the subsidy
adjustment at each decision date during the game and $\Omega ^{f}$ denotes
the set of feasible prices. The government does not give the subsidy to
perpetuity but aims to reach a target of cumulative sales $x_{s}$ with a
minimum expenditure by time ${\tau _{N+1}}<T$, after which it discontinues
the subsidy program. The change in subsidy levels, denoted by $\eta _{i}$,
is done at certain time periods $\tau _{i}$ and the magnitude of change $%
\eta _{i}$ depends on the cumulative sales and subsidy levels such that $%
s(\tau _{i})+\eta _{i}\in \mathcal{S}$. The subsidy levels are constant
between consecutive decision dates $\tau _{i}$ and $\tau _{i+1}$, and the
difference in subsidy levels before and after the intervention time $\tau
_{i}$ is given by: 
\begin{equation}
s\left( \tau _{i}^{+}\right) =s\left( \tau _{i}^{-}\right) +\eta _{i}\quad 
\text{for }i=\{0,1,.....,N\},\;\;\; s(\tau_0)=s_0.
\end{equation}%
The objective of the government is to minimize the expenditure incurred in
reaching the target sales: 
\begin{equation}
\USmodified{J^{g}(p(\cdot ),\eta (\cdot ))}=\min_{\eta _{i}\in \Omega
^{g}(s\USmodified{(\tau_i)}),\,x(\tau _{N+1})\geq x_{s}}\int_{0}^{\tau _{N+1}}e^{-\rho t}s(t)\dot{%
x}(t)dt+\USmodified{\sum_{i=0}^{N+1}e^{-\rho \tau _{i}}C\mathbbm{1}_{\eta _{i}\neq 0}},
\end{equation}
where $\Omega ^{g}(s\USmodified{(\tau_i)})$ denotes the set of feasible subsidy adjustments, with 
$\Omega ^{g}(s):=\{\eta :\eta +\ s\in S\}$ and $C$ is the fixed cost
associated with subsidy adjustments.

To wrap up, we have defined a two-stage differential game model. In the
first stage, the firm and the governments play a noncooperative game,
whereas in the second stage, which starts when the cumulative sales target
is reached, only the firm makes decisions. Consequently, we have to solve an
optimal control problem in the second stage and a differential game in the
first stage. To determine a subgame-perfect equilibrium, we solve the
problem backward. The model involves one state variable and one control
variable for each player. The firm chooses the price of the ZEV in both
stages and the government the subsidy in the first stage. We reiterate that
the firm makes decisions continuously, while the government intervenes only
at some discrete instants of time.
\section{Feedback Stackelberg equilibrium}\label{sec:fse}
\USmodified{We consider a Stackelberg game in which the government, acting as the leader, adopts a feedback strategy $\gamma^g: \{\tau_i\}_{i=1}^N \times \mathbb{R}_+ \times \mathcal{S} \to \Omega^g\bigl(s(\tau_i)\bigr)$
to determine its subsidy adjustment $\eta_i$ at time $\tau_i$, based on the current cumulative sales $x(\tau_i) \in \mathbb{R}_+$ and the prevailing subsidy $s(\tau_i)\in \mathcal S$. We denote by $\mathbb{R}_+$ the set of nonnegative real numbers (the domain for cumulative sales), while $\mathcal{S}$ is the set of all possible subsidy levels, and $\Omega^g\bigl(s(\tau_i)\bigr)$ and $\Omega^f$ are the sets of feasible actions for the government and firm, respectively. The set of all possible government strategies is denoted by $\Gamma^g$. Subsequently, the firm, acting as the follower, observes the government’s choice $\eta_i$ and sets its price according to a feedback strategy $
\gamma^f : [0,T] \times \mathbb{R}_+ \times \mathcal{S} \times \Omega^g \;\to\; \Omega^f,$ 
so that for each $t \in [\tau_i, \tau_{i+1})$, the price is given by 
$p(t) = \gamma^f\bigl(t,\,x(t),\,s(\tau_i),\,\eta_i\bigr).$ 
 The set of all  strategies of the firm is denoted by $\Gamma^f$. By explicitly conditioning the firm’s price choice on $\eta_i$, the model ensures that the follower’s best-response strategy adapts to the leader’s (government’s) subsidy announcements at each decision date $\tau_i$. Furthermore the feedback strategies of both the government ($\gamma^f$)
and the firm ($\gamma^g$) can adapt to the cumulative sales $(x(t))$ and subsidy levels $s(t)$ available at the time of taking an action.}

\begin{definition}
\label{def:fse} We say that $\hat{\gamma}^{f}$ is the firm's best response
to the strategy $\gamma ^{g}$ of the government if \USmodified{%
\begin{align*}
&J^{f}(\hat{\gamma}^{f}(\gamma
^{g}), \gamma^{g}) \geq J^{f}(\gamma^{f}(\gamma
^{g}), \gamma^{g})\;\;\forall \gamma ^{f}\in \Gamma
^{f}.
\end{align*}%
} Similarly, $\hat{\gamma}^{g}$ is the equilibrium strategy of the
government if \USmodified{%
\begin{align*}
&J^{g}(\hat{\gamma}^{f}(\hat{\gamma}^{g}), 
\hat{\gamma}^{g}) \leq J^{g}(\hat{\gamma}^{f}(\gamma^{g}), 
\gamma^{g})\;\;\forall \gamma ^{g}\in \Gamma ^{g}.
\end{align*}%
}
The pair $(\hat{\gamma}^g, \hat{\gamma}^f)$ is called the Feedback
Stackelberg equilibrium (FSE) of the game.
\end{definition}

Once the date to meet the target sales is reached, government stops the
subsidy program. In the \USmodified{next} section, we provide sufficient conditions to
characterize the optimal pricing strategy of the firm \USmodified{during the time interval $[\tau
_{N+1},T]$} after the subsidy program ends.

\subsection{After the end of the subsidy program}

Denote by $v^{f}:[0,T]\times \mathbb{R}_{+}\times \mathcal{S}\rightarrow 
\mathbb{R}$ the value function of the firm.\footnote{\USmodified{The value
function of the firm gives the payoff to go from any position of the game
defined by the couple time and state value, i.e., $( t,x(t), s(t)) $ until
the end of the planning horizon.}} 
After the subsidy program ends, the firm solves a linear-quadratic control
problem, so the value function of the firm satisfies the following
Hamilton-Jacobi-Bellman (HJB) equation: 
\begin{equation*}
\rho v^{f}(t,x(t))-v_{t}^{f}(\USmodified{t},
x(t))=\max_{p(t)}\left[\left(p(t)-c(x(t))+v_{x}^{f}(t,x(t))\right)(\alpha _{1}+\alpha
_{2}x(t)-\beta (p(t)-p_{a}))\right],
\end{equation*}
where $v_{z}^{f}$ is the derivative of $v^{f}$ with respect to variable $z$.
We have suppressed the dependence of $v^f$ on the subsidy since it is
constant between decision dates of the government. Assuming interior
solutions $0<p(t)<\infty $ for $t\in \lbrack 0,T]$, the optimal price
charged by the firm is \USmodified{ obtained using the first-order optimality conditions}, and is given by: 
\begin{equation}
\hat{\gamma}^{f}(t,x(t))=p^*(t)=\frac{1}{2}\left( \frac{\alpha _{1}+\alpha
_{2}x(t)}{\beta }+p_{a}+b_{1}-b_{2}x(t)-v_{x}^{f}(t,x(t))\right) \text{ for }%
t\in \lbrack \tau _{N+1},T].  \label{eq:fse:price_after_subsidy}
\end{equation}

Given the linear-quadratic structure of the problem, we make the informed
guess that the value function is quadratic in the state, that is, $%
v^{f}(t,x(t))=\frac{1}{2}k_{2}(t)x(t)^{2}+k_{1}(t)x(t)+k_{0}(t)$, where $%
k_{j}(t),j=0,1,2$ are the unknown time functions to be determined.
Therefore, the optimal price charged by the firm is given by 
\begin{equation}
p^*(t)=\frac{1}{2}\left( \left( \frac{\alpha _{2}}{\beta }
-b_{2}-k_{2}(t)\right) x(t)+\frac{\alpha _{1}}{\beta }+p_{a}+b_{1}-k_{1}(t)
\right) \text{ for }t\in \lbrack \tau _{N+1},T].
\end{equation}

Substituting the quadratic form of the value function $v^{f}(t,x(t))$ and
optimal price $p^*(t)$ in the HJB equation yields: 
\begin{align*}
\frac{\rho }{2}k_{2}(t)x(t)^{2}&+\rho k_{1}(t)x(t)+\rho k_{0}(t) -\frac{1}{2}\dot{%
k }_{2}(t)x(t)^{2}-\dot{k}_{1}(t)x(t)-\dot{k}_{0}(t) \\
& =\frac{\beta }{4}\left( \frac{\alpha _{1}+\alpha _{2}x(t)}{\beta }
+p_{a}-b_{1}+b_{2}x(t)+k_{2}(t)x(t)+k_{1}(t)\right) ^{2}.
\end{align*}
On comparing the coefficients, we obtain 
\begin{subequations}
\label{eq:ricati_after_subsidy}
\begin{align}
& \rho k_{2}(t)-\dot{k}_{2}(t)=\frac{\beta }{2}\left( \frac{w_{2}}{\beta }
+k_{2}(t)\right) ^{2}, \\
& \rho k_{1}(t)-\dot{k}_{1}(t)=\frac{\beta }{2}\left( \frac{w_{1}}{\beta }
+k_{1}(t)\right) \left( \frac{w_{2}}{\beta }+k_{2}(t)\right), \\
& \rho k_{0}(t)-\dot{k}_{0}(t)=\frac{\beta }{4}\left( \frac{w_{1}}{\beta }
+k_{1}(t)\right) ^{2},
\end{align}
where $w_{1}:=\alpha _{1}+\beta (p_{a}-b_{1})$ and $w_{2}:=\alpha _{2}+\beta
b_{2}$.

\subsection{Before the subsidy ends}

Between the two consecutive subsidy updates at decision dates $\tau _{i}$
and $\tau _{i+1}$, $i=0,1,\ldots ,N$, the value function of the firm evolves
according to the following HJB equation: 
\end{subequations}
\begin{align}
&\rho v^{f}(t,x(t))-v_{t}^{f}(t,x(t))\notag\\&\;=\USmodified{\max_{p(t)\geq
0}}\left[\left(p(t)-c(x(t))+v_{x}^{f}(t,x(t)))(\USmodified{\alpha_1+\alpha_2x(t)} -\beta (p(t)-p_{a}-s(\tau
_{i})\USmodified{-\eta_i})\right)\right].  \label{eq:hjb}
\end{align}%
Assuming interior solutions for the optimal price $p^*(t)$ and
quadratic form of the value function, the optimal price for $t\in (\tau
_{i},\tau _{i+1})$ charged by the firm is obtained using the first-order
condition: 
\begin{equation}
p^*(t)=\frac{1}{2}\left(
\left( \frac{\alpha _{2}}{\beta }-b_{2}-k_{2}(t)\right) x(t)+\frac{\alpha
_{1}}{\beta }+p_{a}+s(\tau _{i}) \USmodified{ + \eta_i} +b_{1}-k_{1}(t)\right) .
\label{eq:fse:price}
\end{equation}%
Substituting $p^*(t)$ in the state dynamics \eqref{eq:state_dyn}, the
state evolution over $(\tau _{i},\tau _{i+1})$ is given by 
\begin{equation}
\dot{x}(t)=\frac{1}{2}\left( \frac{w_{1}}{\beta }+\beta (s(\tau
_{i}) + \USmodified{\eta_i}-k_{1}(t))+\left( \frac{w_{2}}{\beta }+\beta k_{2}(t)\right)
x(t)\right)  \label{state:eq_opt}
\end{equation}%
Substituting the quadratic form of the value function $v^{f}(t,x(t))$ and
optimal price $p^*(t)$ in the HJB equation, and comparing the
coefficients yields for $i\in \{0,1,\ldots ,N\}$: 
\begin{subequations}
\label{eq:riccati}
\begin{align}
& \rho k_{2}(t)-\dot{k}_{2}(t)=\frac{\beta }{2}\left(k_{2}(t)+\frac{w_{2}}{\beta }\right)^{2}\text{ for }t\in (\tau _{i},\tau _{i+1}), \\
& \rho k_{1}(t)-\dot{k}_{1}(t)=\frac{\beta }{2}\left( \frac{w_{1}}{\beta }%
+k_{1}(t)+s(\tau _{i})+\USmodified{\eta_i}\right) \left( \frac{w_{2}}{\beta }+k_{2}%
\USmodified{(t)}\right) \text{ for }t\in (\tau _{i},\tau _{i+1}), \\
& \rho k_{0}(t)-\dot{k}_{0}(t)=\frac{\beta }{4}\left( \frac{w_{1}}{\beta }%
+k_{1}(t)+s(\tau _{i})+\USmodified{\eta_i}\right) ^{2}\text{ for }t\in (\tau _{i},\tau _{i+1}). \label{eq:k0}
\end{align}%
The continuity of the value function at the time instants $\tau
_{i},i=\{0,1,\ldots ,N+1\}$ yields the following relation: 
\end{subequations}
\begin{equation*}
\frac{1}{2}k_{2}(\tau _{i}^{-})x(\tau _{i}^{-})^{2}+k_{1}(\tau
_{i}^{-})x(\tau _{i}^{-})+k_{0}(\tau _{i}^{-})=\frac{1}{2}k_{2}(\tau
_{i}^{+})x(\tau _{i}^{+})^{2}+k_{1}(\tau _{i}^{+})x(\tau
_{i}^{+})+k_{0}(\tau _{i}^{+}).
\end{equation*}%
From the continuity of $x(\cdot )$ in $t$, we obtain: 
\begin{equation}
k_{m}(\tau _{i}^{-})=k_{m}(\tau _{i}^{+}),\text{ for }m=\{0,1,2\}. \label{eq:kink:km}
\end{equation}%
Therefore, the value function has a kink at the decision date $\tau _{i}$ if
the subsidy adjustment is made. \USmodified{Since $k_{2}(t)$ does not depend on the
subsidy levels, we can analytically characterize its solution (see Appendix %
\ref{appen_k2}): %
When $\rho ^{2}-2\rho \,w_{2}\geq 0,$ we have 
\begin{subequations}
\label{eq:k2}
\begin{equation}
k_{2}(t)=\frac{\rho -w_{2}}{\beta }+\frac{\Delta }{\beta }\,\tanh \!\left[ 
\frac{\Delta }{2}\left( t-T+\frac{2}{\Delta }\operatorname{arctanh}\!\left( \frac{%
w_{2}-\rho }{\Delta }\right) \right) \right], \forall t\in [0,T],
\end{equation}%
where $\Delta =\sqrt{\rho ^{2}-2\rho \,w_{2}}.$ When $2\rho \,w_{2}-\rho
^{2}\geq 0$, we get 
\begin{equation}
k_{2}(t)=\frac{\rho -w_{2}}{\beta }+\frac{\Theta }{\beta }\tan \left[ -\frac{%
\Theta }{2}\left( t-T-\frac{2}{\Theta }\operatorname{arctan}\left( \frac{w_{2}-\rho 
}{\Theta }\right) \right) \right], \forall t\in [0,T],
\end{equation}%
where $\Theta =\sqrt{2\rho \,w_{2}-\rho ^{2}}$. Similarly, we can also
derive a recurrence relationship for $k_{1}(t)$ (see Appendix \ref{appen_k1}%
). When $\rho ^{2}-2\rho \,w_{2}\geq 0$:
\end{subequations}
\begin{subequations}
\label{eq:k1}
\begin{align}
k_{1}(t)&=\frac{1}{I_{\Delta }(t)}\left[ I_{\Delta }(\tau_{i+1})k_{1}(\tau _{i+1})+\rho \left( s(\tau _{i})+\eta_i+\frac{w_{1}}{%
\beta }\right) \int_{t}^{\tau _{i+1}}I_{\Delta }(u)\,du\right]\notag\\&\hspace{1cm} +(I_{\Delta }(\tau_{i+1})/I_{\Delta}(t)-1)(s(\tau_{i})+\eta_i+w_{1}/\beta),\, \forall t\in (\tau_i,\tau_{i+1}),
\end{align}%
where 
\begin{equation*}
I_{\Delta }(t)=\exp \!\left( -\frac{\rho }{2}\,t\right) \cosh \!\left( \frac{%
\Delta }{2}\left( t-T+\frac{2}{\Delta }\operatorname{arctanh}\!\left( \frac{w_{2}-\rho 
}{\Delta }\right)\right) \right) ,
\end{equation*}%
while for $2\rho \,w_{2}-\rho ^{2}\geq 0$: 
\begin{align}
k_{1}(t)=&\frac{1}{I_{\Theta }(t)}\left[I_{\Theta }(\tau_{i+1}) k_{1}(\tau _{i+1})+\rho \left( s(\tau _{i})+\eta_i+\frac{w_{1}}{%
\beta }\right) \int_{t}^{\tau _{i+1}}I_{\Theta }(u)\,du\right]\notag \\&\hspace{1cm}+(I_{\Theta }(\tau_{i+1})/I_\Theta(t)-1)(s(\tau _{i})+\eta_i+w_{1}/\beta),\, \forall t\in (\tau_i,\tau_{i+1}),  \label{eq:k1:recu}
\end{align}%
where 
\begin{equation*}
I_{\Theta }(t)=\exp \!\left( -\frac{\rho }{2}t\right) \Bigg|\cos \!\left( 
\frac{\Theta }{2}\left(t-T-\frac{2}{\Theta }\operatorname{arctan}\left( \frac{w_{2}-\rho 
}{\Theta }\right) \right)\right) \Bigg|.
\end{equation*}
\end{subequations}
Note that we can obtain $k_{1}(\tau _{i})$ using the recurrence
equation: when $\rho ^{2}-2\rho w_{2}\geq 0,$ we have 
\begin{subequations}
\begin{align}
k_{1}(\tau _{i})=&\frac{1}{I_{\Delta }(\tau _{i})}\left[ I_{\Delta }(\tau _{i+1})k_{1}(\tau
_{i+1})+\rho \left( s(\tau
_{i})+\eta_i+w_{1}/\beta\right) \int_{\tau _{i}}^{\tau
_{i+1}}I_{\Delta }(u)\,du\right] \notag\\&+(I_{\Delta }(\tau_{i+1})/I_\Delta(\tau_i)-1)(s(\tau _{i})+\eta_i+w_{1}/\beta),
\label{eq:k1tau_recu_del}
\end{align}%
and when $2\rho w_{2}-\rho ^{2}\geq 0,$we obtain 
\begin{align}
k_{1}(\tau _{i})=&\frac{1}{I_{\Theta }(\tau _{i})}\left[ I_{\Theta }(\tau _{i+1})k_{1}(\tau
_{i+1})+\rho \left( s(\tau
_{i})+\eta_i+w_{1}/\beta\right) \int_{\tau _{i}}^{\tau
_{i+1}}I_{\Theta }(u)\,du\right]\notag\\&+(I_{\Theta }(\tau_{i+1})/I_\Theta(\tau_i)-1)(s(\tau _{i})+\eta_i+w_{1}/\beta),
\label{eq:k1tau_recu_theta}
\end{align}
\end{subequations}
where we use the fact that $k_1(T)=0.$
On substituting the expressions for $k_1(t)$ and $k_2(t)$ in the expression of price given in  \eqref{eq:fse:price}, we arrive at the following result:
\begin{proposition}\label{eq:prop}
The price charged by the firm for $t\in (\tau_i, \tau_{i+1}), i\leq N$ is
as follows: when $\rho ^{2}-2\rho w_{2}\geq 0$, we have 
\begin{subequations}
\begin{align}
&\hat{\gamma}^{f}(t, x(t),s(\tau _{i}),\eta_i) \notag\\&=\frac{1}{2\beta }\left( -\rho -\Delta \tanh \!\left[ \frac{\Delta }{2}\left( t-T+\frac{2}{\Delta }%
\operatorname{arctanh}\!\left( \frac{w_{2}-\rho}{\Delta }\right)\right) \right]
\right) x(t)  \notag \\
& \hspace{1cm}-\frac{1}{2I_{\Delta }(t)}\left[ I_{\Delta }(\tau_{i+1})k_{1}(\tau _{i+1})+\rho \left( s(\tau _{i})+\eta_i+\frac{w_{1}}{%
\beta }\right) \int_{t}^{\tau _{i+1}}I_{\Delta }(u)\,du\right]\notag\\&\hspace{2cm} +(2-I_{\Delta }(\tau_{i+1})/I_{\Delta}(t))(s(\tau_{i})+\eta_i+w_{1}/\beta)/2,
\label{p1_less}
\end{align}%
and when $2\rho w_{2}-\rho ^{2}\geq 0$, we have:
\begin{align}
&\hat{\gamma}^{f}(t, x(t),s(\tau _{i}), \eta_i)\notag\\& =\frac{1}{2\beta }\left( -\rho -\Theta \tan \!\left[ -\frac{\Theta }{2}\left( t-T-\frac{2}{%
\Theta }\operatorname{arctan}\!\left( \frac{w_{2}-\rho }{\Theta }\right) \right) %
\right] \right) x(t)  \notag \\
& \hspace{1cm}-\frac{1}{2I_{\Theta }(t)}\left[I_{\Theta }(\tau_{i+1}) k_{1}(\tau _{i+1})+\rho \left( s(\tau _{i})+\eta_i+\frac{w_{1}}{%
\beta }\right) \int_{t}^{\tau _{i+1}}I_{\Theta }(u)\,du\right]\notag \\&\hspace{2cm}+(2-I_{\Theta }(\tau_{i+1})/I_\Theta(t))(s(\tau _{i})+\eta_i+w_{1}/\beta)/2, 
\label{p1_high}
\end{align}%
where $w_{1}:=\alpha _{1}+\beta (p_{a}-b_{1})$ and $w_{2}:=\alpha _{2}+\beta
b_{2}$.

\noindent
The price for $t\in \lbrack \tau _{N+1},T]$ when $\rho ^{2}-2\rho w_{2}\geq 0$ is given by: 
\end{subequations}
\begin{subequations}
\begin{align}
\hat{\gamma}^{f}(t,x(t)) &=\frac{1}{2\beta }\left( -\rho
-\Delta \tanh \!\left[ \frac{\Delta }{2}\left( t-T+\frac{2}{\Delta }\operatorname{%
arctanh}\!\left( \frac{w_{2}-\rho}{\Delta }\right)\right) \right]
\right) x(t)  \notag \\
& \;\;-\frac{\rho w_1}{2\beta I_{\Delta }(t)}\int_{t}^{T}I_{\Delta }(u)\,du +\frac{(2-I_{\Delta }(T))w_{1}}{2\beta I_\Delta(t)},
\label{p1_low_nosubs}
\end{align}
and when $2\rho w_{2}-\rho ^{2}\geq 0$, the price for $t\in \lbrack \tau _{N+1},T]$ is given by:
\begin{align}
\hat{\gamma}^{f}(t,x(t)) &=\frac{1}{2\beta }\left( -\rho
-\Theta \tan \!\left[ -\frac{\Theta }{2}\left( t-T-\frac{2}{\Theta }\operatorname{%
arctan}\!\left( \frac{w_{2}-\rho }{\Theta }\right) \right) \right] \right)
x(t)  \notag \\
& \;\;-\frac{\rho w_1}{2\beta I_{\Theta }(t)}\int_{t}^{T}I_{\Theta }(u)\,du +\frac{(2-I_{\Theta }(T))w_{1}}{2\beta I_\Theta(t)}.
\label{p1_high_nosubs}
\end{align}%
\end{subequations}
\end{proposition}

Since $\operatorname{cosh}(y)>0$ for all $y$, we obtain $I_{\Delta }(t)\geq
0 $, and clearly $I_{\Theta }(t)\geq 0\,\forall t$, which implies that the
coefficient of $s(\tau _{i})+\eta_i$ is always positive. Therefore, the price
charged by the firm is increasing linearly with the subsidy levels and is
linear in the cumulative sales. The result that the price is increasing in
the subsidy is in line with what has been obtained in the empirical
literature that assessed the impact of vehicle scrappage programs and cited
in the introduction. }

Next, we consider the control problem of the government. Let the value
function of the government be denoted by $v^{g}:[0,T]\times \mathbb{R}%
_{+}\times \mathcal{S}\rightarrow \mathbb{R}$. Government would stop the
subsidy program at $\tau _{N+1}$. \USmodified{We define the continuation set, representing the scenario where the government has met its target by time $\tau_{N+1}$:}
\begin{equation*}
\mathcal{C}=\{(\tau _{N+1},x(\tau _{N+1}), \USmodified{s(\tau_{N+1})}):x(\tau _{N+1})\geq x_{s}\}.
\end{equation*}%

Equivalently, the value function of the government satisfies the following
 condition at time $\tau _{N+1}:$ 
\begin{equation}
v^{g}(\tau _{N+1}, x(\tau _{N+1}), s(\tau _{N+1}))=%
\begin{cases}
0,\;\text{ if }\;(\tau _{N+1},x(\tau _{N+1}), \USmodified{s(\tau _{N+1})})\in \mathcal{C}, \\ 
\infty ,\;\text{ otherwise. }%
\end{cases}
\label{eq:terminal_cond}
\end{equation}%
\USmodified{We denote by $\mathcal{M}v^{g}$ the minimum cost-to-go (evaluated at time $\tau_i$)  for adjusting the subsidy so as to reach
the target $x_{s}$ by time $\tau _{N+1}$} 
and is defined by: 
\begin{align}
& \mathcal{M}v^{g}(\tau _{i},x(\tau _{i}), s(\tau _{i}))\notag\\&=\min_{\eta _{i}\in \Omega
^{g}(s(\tau _{i}))}\Big(\int_{\tau _{i}^{+}}^{\tau _{i+1}^{-}}\USmodified{e^{-\rho (t-\tau_i)}}\left( s(\tau _{i})+\eta _{i}\right)
\dot{x}(t)dt+C\USmodified{\mathbbm{1}_{\eta _{i}\neq 0}}+\USmodified{e^{-\rho (\tau_{i+1}-\tau_{i})}}v^{g}(\tau _{i+1},x(\tau
_{i+1}),s\USmodified{(\tau_{i+1})})\Big)  \notag \\
& =\min_{\eta _{i}\in \Omega ^{g}(s(\tau _{i}))}\Bigg(\int_{\tau _{i}^{+}}^{\tau
_{i+1}^{-}}\USmodified{e^{-\rho (t-\tau_i)}}\frac{s(\tau _{i})+\eta _{i}}{2}\left( \frac{w_{1}}{\beta }+\beta
(s(\tau _{i})+\eta _{i}-k_{1}(t))+\left( \frac{w_{2}}{\beta }+\beta k_{2}(t)\right)
x(t)\right) dt  \notag \\
& \hspace{1cm}+C\USmodified{\mathbbm{1}_{\eta _{i}\neq 0}}+\USmodified{e^{-\rho (\tau_{i+1}-\tau_{i})}}v^{g}(\tau _{i+1},x(\tau
_{i+1}),s\USmodified{(\tau_{i+1})}\Bigg),  \label{eq:op_M}
\end{align}%
where $\mathcal{M}$ is the intervention operator, \USmodified{and the last equality is obtained by using \eqref{state:eq_opt}.} By construction, the
value function satisfies the following relation at each decision date, $\tau_i, \, i=\{0,1, \ldots, N+1\}$: 
\begin{equation}
v^{g}(\tau _{i},x(\tau _{i}), s(\tau _{i}))=\mathcal{M}%
v^{g}(\tau _{i},x(\tau _{i}), s(\tau _{i})).
\label{eq:optim_imp}
\end{equation}

\begin{theorem}
Suppose there exists a  function $v^{f}=\frac{1}{2}%
k_{2}\USmodified{(t)}x(t)^{2}+k_{1}\USmodified{(t)}x(t)+k_{0}\USmodified{(t)}$ such that $k_{m}:[0,T]\rightarrow \mathbb{R}$
for $m\in \{0,1,2\}$ satisfy \eqref{eq:ricati_after_subsidy}, %
\eqref{eq:riccati}, and \eqref{eq:kink:km} for $t\in (\tau _{j},\tau _{j+1})$%
, $j\in \{0, 1, \ldots ,N+1\}$,  $\tau _{0}:=0$ \USmodified{and $\tau_{N+2}:=T$}. Furthermore, suppose there a
exists a function $v^{g}$ that satisfies \eqref{eq:terminal_cond}, %
\eqref{eq:op_M}, and \eqref{eq:optim_imp}. Then, $\hat{\gamma}^{f}(\cdot )$
\USmodified{given in Proposition \ref{eq:prop}} and the
subsidy adjustments $\hat{\gamma}^{g}(\cdot )$ defined below constitute
the FSE strategy of the firm and government, respectively, \USmodified{and $v^f$ and $v^g$ are their respective value functions}: 
\begin{align}
&\hat{\gamma}^{g}(t, x(t), s(\tau_{j}))=\arg \min_{\eta \in \Omega ^{g}(s(\tau_{j}))}\Bigg[%
\int_{\tau _{j}^{+}}^{\tau _{j+1}^{-}}\Big(e^{-\rho (t-\tau_j)}\frac{s(\tau _{j})+\eta }{2}%
\Big( \frac{w_{1}}{\beta }+\beta (s(\tau _{j})+\eta
-k_{1}(t))\notag\\&\hspace{0.5cm}+\left( \frac{w_{2}}{\beta }+\beta k_{2}(t)\right) x(t)\Big)\Big) dt +C\USmodified{\mathbbm{1}_{\eta _{j}\neq 0}}+\USmodified{e^{-\rho (\tau_{j+1}-\tau_{j})}}v^{g}(\tau _{j+1},x(\tau
_{j+1}),s(\tau _{j+1})) \Bigg].  \label{eq:fse:subsidy} 
\end{align}
\end{theorem}
\begin{proof}
From Definition \ref{def:fse}, we will show that for all $j\in \{0, 1, \ldots, N+1\}$, we have:
\USmodified{\begin{align*}
& v^{f}(\tau _{j}, x(\tau_j), s(\tau_j))=J^{f}(
\hat{\gamma}^{f}(\cdot, 
\hat{\gamma}^{g}),  \hat{\gamma}^{g}(\cdot)),
\\&v^{g}(\tau _{j},x(\tau_j),s(\tau_j))=J^{g}(
\hat{\gamma}^{f}(\cdot, 
\hat{\gamma}^{g}),  \hat{\gamma}^{g}(\cdot)), \\
& v^{f}(\tau _{j},x(\tau_j), s(\tau_j) )\geq J^{f}(
\gamma^{f}(\cdot, 
\hat{\gamma}^{g}),  \hat{\gamma}^{g}(\cdot)),  \,\forall \gamma ^{f}\in \Gamma ^{f}, \\
& v^{g}(\tau _{j},x(\tau_j),s(\tau_j))\leq J^{g}(\hat{\gamma}^{f}(\cdot,
\gamma^{g}),  \gamma^{g}(\cdot)),  \,\forall \gamma ^{g}\in \Gamma ^{g}.
\end{align*}}
Suppose the feedback strategy of the
government is $\hat{\gamma}^{g}$ so that the subsidy change
at decision date $\tau _{i}$ is $\eta _{i}=\gamma^{g}(\tau _{i},x(\tau _{i}), s(\tau _{i})),\, \forall i\geq j$. \USmodified{Let $\gamma^{f}$ be an arbitrary feedback strategy of the
firm  specifying that at any time $t$, the price  is given by $p(t)=\gamma^{f}(t, x(t), s(\tau_i), \eta_i)$, where $\tau_i$ is the earliest decision date that is greater than or equal to $t$.}  Let the corresponding
state trajectory be denoted by $x_{1}(\cdot )$. Using the total derivative
of $e^{-\rho t}v^{f}(\cdot )$ between $(\tau _{i-1},\tau _{i})$, integrating
with respect to $t$ from $\tau _{i-1}$ to $\tau _{i}$, and taking the
summation for all $i\geq \USmodified{j}$, we obtain 
\begin{align}
& e^{-\rho T}v^{f}(T,x_{1}(T), \USmodified{s(T)})-e^{-\rho \tau _{j}}v^{f}(\tau _{j},x(\tau
_{j}), s(\tau_j))  \notag \\
& =\sum_{i\,\geq\, j}\int_{\tau _{i}^{+}}^{\tau _{i+1}^{-}}e^{-\rho h}\Big(
-\rho
v^{f}(h,x_{1}(h), s(\tau_i)+\eta_i)+v_{h}^{f}(h,x_{1}(h), s(\tau_i)+\eta_i)\notag\\&\hspace{1cm}+v_{x}^{f}(h,x_{1}(h), s(\tau_i)+\eta_i)\USmodified{\dot{x}_1(h)}\Big) dh,
\label{eq:total_diff}
\end{align}%
\USmodified{where $\dot{x}_1(h) = \alpha _{1}+\alpha _{2}x_1(h)-\beta (p(h)-s(\tau_i)-\eta_i-p_{a})$ for $h\in (\tau_i, \tau_{i+1})$.}
From the HJB equation \eqref{eq:hjb}, the following inequality holds for $%
(\tau _{i},\tau _{i+1})$: 
\begin{align*}
&-\rho
v^{f}(h,x_{1}(h), s(\tau_i)+\eta_i)+v_{h}^{f}(h,x_{1}(h), s(\tau_i)+\eta_i)\notag\\&\hspace{1cm}+v_{x}^{f}(h,x_{1}(h), s(\tau_i)+\eta_i)\USmodified{\dot{x}_1(h)}\leq -(p(h)-c(x_{1}(h)))\dot{x}_1(h).
\end{align*}%
Using the above inequality in \eqref{eq:total_diff} and substituting $%
v^{f}(T,x_{1}(T), s(T))=0$, we obtain: 
\begin{equation*}
-e^{-\rho \tau _{j}}v^{f}(\tau _{j},x(\tau _{j}), s(\tau_j))\leq -\sum_{i\geq
j}\int_{\tau _{i}^{+}}^{\tau _{i+1}^{-}}e^{-\rho h}(p(h)-c(x_{1}(h)))\dot{x}_1(h)dh.
\end{equation*}%
Rearranging the above equation yields: 
\begin{align*}
v^{f}(\tau _{j},x(\tau _{j}), s(\tau_j))&\geq \sum_{i\geq j}\int_{\tau _{i}^{+}}^{\tau
_{i+1}^{-}}e^{-\rho (h-\tau _{j})}(p(h)-c(x_{1}(h))\dot{x}_1(h)dh\\&\;\;=J^{f}\left(\gamma^{f}(\cdot, \hat{\gamma}^{g}),  \hat{\gamma}^{g}\right).
\end{align*}%
\USmodified{The inequality becomes an equality when the price is set as specified in Proposition \ref{eq:prop}, meaning that $\hat{\gamma}%
^{f}(\cdot, \hat{\gamma}^g)$ constitutes the firm’s best response strategy to the government’s subsidy adjustment strategy $\hat{\gamma}^{g}$.}

Next, we consider the intervention problem of the government. \USmodified{For a feedback strategy $\hat{\gamma}^f$ of the firm,  let $\gamma^{g}$ denote the arbitrary feedback
strategy  of the government, which adjusts the subsidy by an amount $a_j$ at each decision date $\tau_j$.}  Let $x_{2}(t)$ denote the
state corresponding to the cumulative sales and $s(\tau_j)$ denotes the subsidy at any decision date $\tau_j$. Since \eqref{eq:op_M} and %
\eqref{eq:optim_imp} hold for all $i\geq j$, we obtain:
\begin{align*}
v^g(\tau _{i},x(\tau _{i}), s(\tau _{i}))  &\leq \frac{s(\tau _{i})+\USmodified{a _{i}}}{2}\int_{\tau
_{i}^{+}}^{\tau _{i+1}^{-}}\USmodified{e^{-\rho (t-\tau_i)}}\Big( \frac{w_{1}}{\beta }+\beta
(s(\tau_i)-k_{1}(t)+\USmodified{a _{i}})\\&+\left( \frac{w_{2}}{\beta }+\beta
k_{2}(t)\right) x_{2}(t)\Big) dt +C\mathbbm{1}_{a _{i}\neq 0}+\USmodified{e^{-\rho (\tau_{i+1}-\tau_i)}}v^g(\tau
_{i+1}, x(\tau _{i+1}),s(\tau_{i+1})).
\end{align*}
Since the above relation holds for all $i\in \{j,j+1, \dots, N\}$ \USmodified{and arbitrary adjustments $a _{j}$}, we obtain:
\[ v^g(\tau _{j},x(\tau _{j}), s(\tau _{j}))\leq J^{g}\left(\hat{%
\gamma}^{f}(\cdot, \gamma^g),\gamma^{g}\right)\; \forall \gamma^g \in \Gamma^g.\]

The above inequality is satisfied with an equality for the equilibrium
strategy \eqref{eq:fse:subsidy} of the government. \USmodified{Therefore, $\hat{\gamma}^g$ and $\hat{\gamma}^f$ are best response strategies of the government and the firm, with corresponding value functions $v^g$ and $v^f$.}
\end{proof}

\USmodified{The above theorem provides a procedure to determine the equilibrium
subsidy plan of the government at the decision dates $\tau_{i}$ for any
given cumulative sales $x(\tau_j)$ and subsidy levels $s(\tau_j)$, which may not lie on the equilibrium path.} 
\begin{algorithm}
\caption{Computing the FSE subsidies}
\begin{algorithmic}[1]
\State  \USmodified{Compute $k_2(t)$ and $k_1(t)$ using \eqref{eq:k2} and \eqref{eq:k1} for $t\in [\tau_{N+1}, T]$, respectively, by setting $s(\tau_{N+1})=0$ and $\eta_{N+1}=0$}
\State Discretize the state $x$ over $\mathcal{G} \in [x_0, x_M]$
\State Set \USmodified{$v^g(\tau_{N+1},  x(\tau_{N+1}), s(\tau_{N+1}))$} according to \eqref{eq:terminal_cond}
\For{$j=N$ to $0$}
    \For{each state $x(\tau_j)\in \mathcal{G}$ at time $\tau_j$}
       \For{each state $s\in \mathcal{S}$ at time $\tau_j$} 
        \State \USmodified{Compute $k_2(t)$ and $k_1(t)$ using \eqref{eq:k2} and \eqref{eq:k1}  for $t\in [\tau_{j}, \tau_{j+1}]$}
            \State Use \eqref{state:eq_opt} to solve for $x(t)$  in $t\in [\tau_{j}, \tau_{j+1})$
            \State Solve \eqref{eq:op_M} to obtain $\mathcal{M}v^{g}(\tau _{j},x(\tau _{j}), s(\tau_j))$ and corresponding  $\eta_j$. 
            \State Set $v^{g}(\tau _{j},x(\tau _{j}), s(\tau_j))\gets \mathcal{M}v^{g}(\tau _{j},x(\tau _{j}),s(\tau_j))$ and $\gamma^g(\tau_j, x(\tau_j), s(\tau_j)) \gets \eta_j$
            \EndFor
    \EndFor
    \State Approximate $v^g(\tau_j,x(\tau_j), s(\tau_j))$ for each $x(\tau_j)\in[ x_0, x_M]$ using interpolation.
\EndFor
\State Return $v^g(0, x_0, s_0)$ for the initial state $x_0$ and $s_0$.
\end{algorithmic}
\label{alg}
\end{algorithm}

To compute the FSE, we can obtain the equilibrium subsidies by dynamic
programming using Algorithm \ref{alg}, \USmodified{where we assume that the cumulative sales lie in the interval $[x_0, x_M]$ for $t\in [0, T]$}. Then, we can use the equilibrium
subsidies to compute \USmodified{the equilibrium price in Proposition \ref{eq:prop}}.

\USmodified{In the next section, we develop a numerical example to illustrate
the kind of results and insights that can be obtained when using our model.
Also, we run a sensitivity analysis to assess the impact of main parameter
values on the players' strategies, that is, the firm's price and the
government's subsidy.} We will exactly solve the models by
enumerating all combinations of subsidy levels at each date, and solving the
optimal-control problem of the firm to obtain $k_{2}(t)$ and $k_{1}(t)$ using \eqref{eq:k2} and \eqref{eq:k1}, respectively. Then, we compute the government's cost in each case by substituting the resulting sales into the government's objective function, and choose the optimal subsidy levels for which the total cost over $%
[0,\tau _{N+1}]$ is minimized where the cost of failing to meet the target sales by 
$\tau _{N+1}$ is taken to be infinite.

\section{\USmodified{Numerical example}}

\label{sec:numerical} Even with the simplest possible specifications of the
functions involved in the model, one cannot solve analytically an iDSG.

As a benchmark case, we adopt the following parameter values:\footnote{%
We ran a large number of numerical examples and the results are
(qualitatively) robust to what we present here.}%
\begin{eqnarray*}
\text{Demand parameters}\text{: } &&\alpha _{1}=6,\text{ \ }\alpha _{2}=0.01,%
\text{ \ }\beta ={0.12},\text{ \ }p_{a}=1,\text{ \ }x_{0}=%
{15}, \\
\text{Cost parameters}\text{: } &&b_{1}={55},\text{ \ \ }%
b_{2}=0.8,\text{ \ }C=10, \\
\text{Other parameters}\text{: } &&T={15},\text{ \ \ }\rho =0.1,%
\text{ \ \ }x_{s}=40.\text{ }
\end{eqnarray*}%
Let the feasible subsidy set be given by $\mathcal{S}=\{0,5,10, 15\}$ and
the subsidy adjustment be made at instants of time $\tau _{1}=0$ and $\tau
_{2}=5$. The subsidy program stops at $\tau _{3}=10$, which is different
from the firm's planning horizon, set here to $T={15}$. 
Based on the results, we answer here our research questions \ref{q1}, \ref%
{q2}, and \ref{q3}. 
\begin{figure}[]
\centering
\begin{subfigure}{0.3\linewidth}
            \centering
\includegraphics[width=\textwidth]{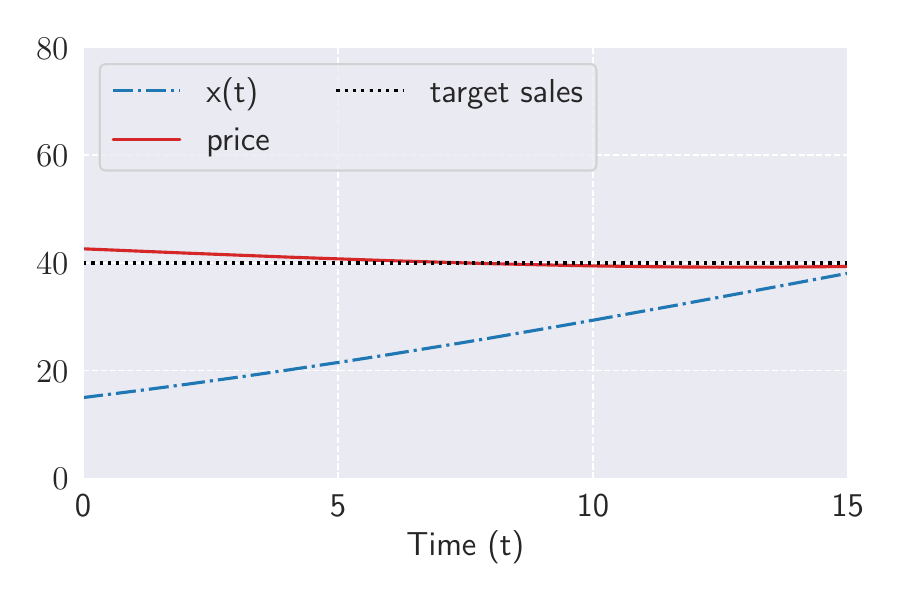}
            \caption{Without subsidy}
            \label{fig:no_subsidy}
        \end{subfigure} 
\begin{subfigure}{0.3\linewidth}
            \centering
\includegraphics[width=\textwidth]{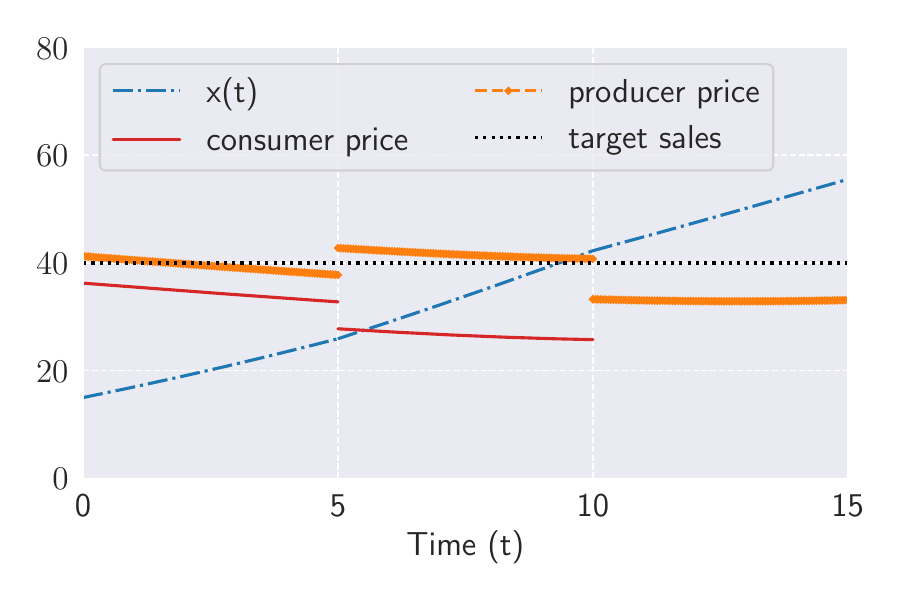}
            \caption{With subsidy}
            \label{fig:price_subsidy}
        \end{subfigure} %
\begin{subfigure}{0.3\linewidth}
            \centering
            \includegraphics[width=\textwidth]
          {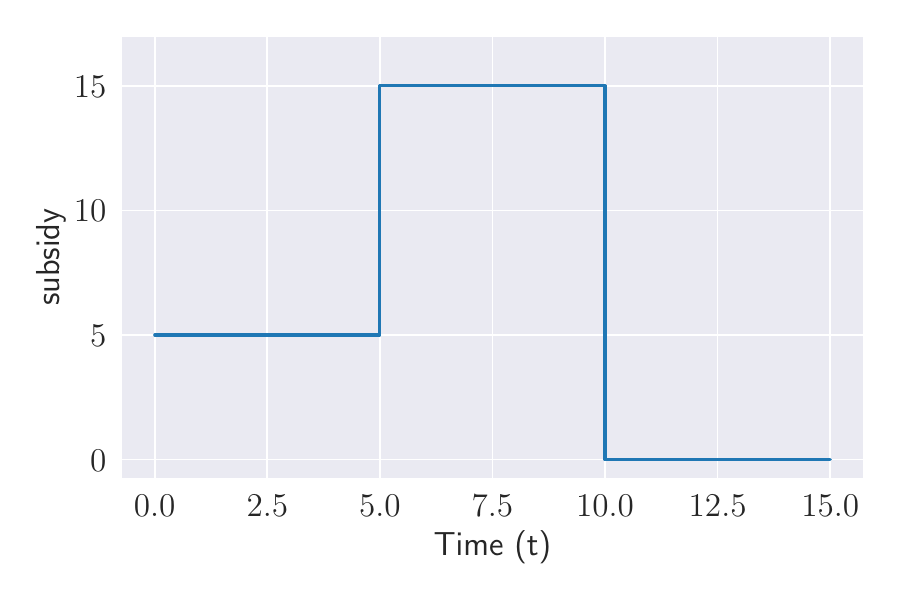}
            \caption{Subsidy}
            \label{fig:subsidylevel}
        \end{subfigure}
\caption{Cumulative sales, equilibrium producer and consumer price, and
equilibrium subsidy for the benchmark case.}
\label{fig:prices}
\end{figure}

\textbf{Price and subsidy}. Figures \ref{fig:no_subsidy} and \ref%
{fig:price_subsidy} show that the consumer's price during the subsidy period
is significantly lower than what she would have paid without the subsidy.
After the subsidy, the price difference is, however, only slightly lower in
the subsidy scenario than in the case without subsidy. In both scenarios,
the seller's price is decreasing over time, which is the consequence of
learning-by-doing and word-of-mouth effects.

The government sets the subsidy level at $5$ at time $0$ and at $15$
at time $5$ (see Figure \ref{fig:subsidylevel}). Starting with a high
subsidy, slightly more than ${10}\%$ of the seller's price, is
meant to trigger a snowballing effect in the adoption process. Indeed, high
subsidy leads to high demand, which accelerates the reduction in the unit
production cost and the positive word-of-mouth effect. In turn, the price
goes down and adoption rate up.

Comparing the results in Figures \ref{fig:no_subsidy} and \ref%
{fig:price_subsidy} shows that the adoption rate is higher when the
government offers a subsidy than when it does not, which is expected. Note
that the target would not been reached without subsidy. 
\paragraph{Firm's strategic behavior.}
The subsidy period is given by the time interval $I=\left[ 0,{10}%
\right].$ Denote by $p^{s}(t)$ the firm's price when the
government offers a subsidy and by $p(t)$ the price when it
does not. Our results show that $p^{s}(t) >p(t)$
for all $t\in [5, 10]$. Then, the answer to our second research question is
unambiguous: the firm indeed takes advantage of the implementation of the
subsidy program to increase its price during the time interval $[5, 10]$. One
implication of this strategic behavior is that the subsidy program is not
achieving its full potential in terms of decreasing the price and raising
the adoption rate. The same conclusion was reached in \cite%
{kaul2016incidence} and \cite{Jimenez2016} in their evaluation of
the car scrappage program.

\paragraph{Cost and benefits.}

Whereas the cost of a subsidy program is straightforward to compute, its
benefits are more complex to determine. Here, the subsidy program costs
taxpayers $1593.29$. The benefits can be assessed in terms of
consumer surplus, firm's profit, and environmental impact.
\begin{figure}[H]
\centering
\begin{subfigure}{0.35\linewidth}
            \centering
\includegraphics[width=\textwidth]{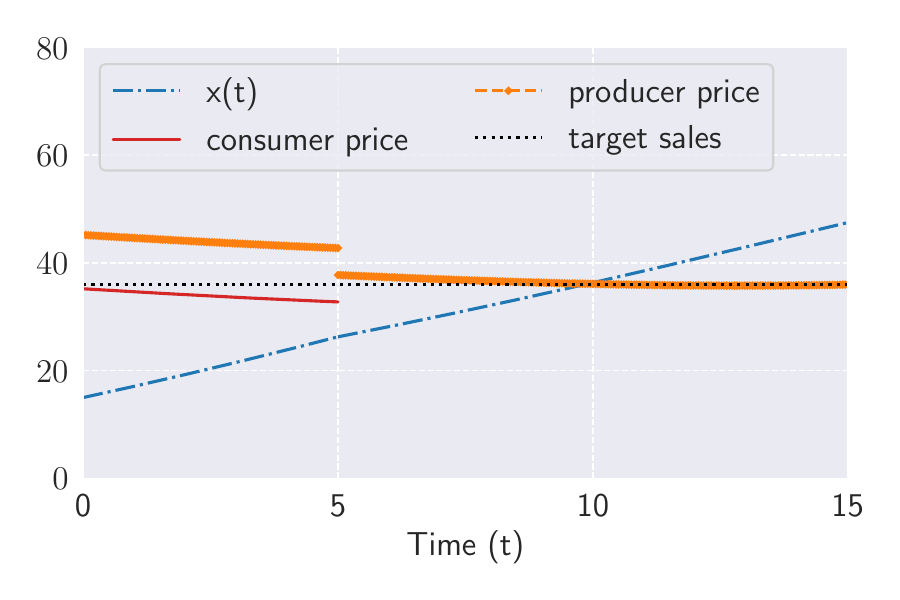}
            \caption{Target = 36}    \label{fig:price_target_lowest}
        \end{subfigure} \quad 
\begin{subfigure}{0.35\linewidth}
            \centering
\includegraphics[width=\textwidth]{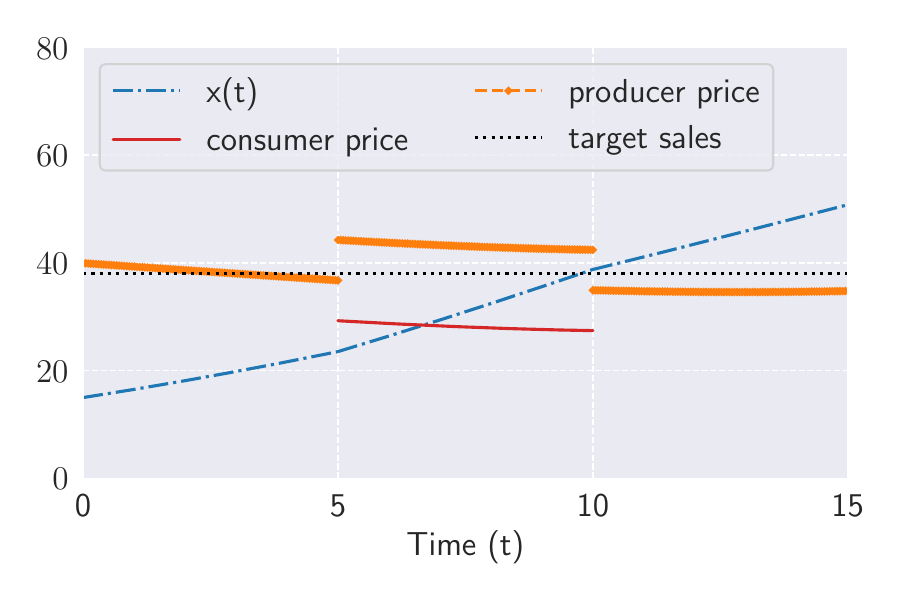}
            \caption{Target = 38}
            \label{fig:price_subsidy_lower}
        \end{subfigure} \\
\begin{subfigure}{0.35\linewidth}
            \centering
\includegraphics[width=\textwidth]{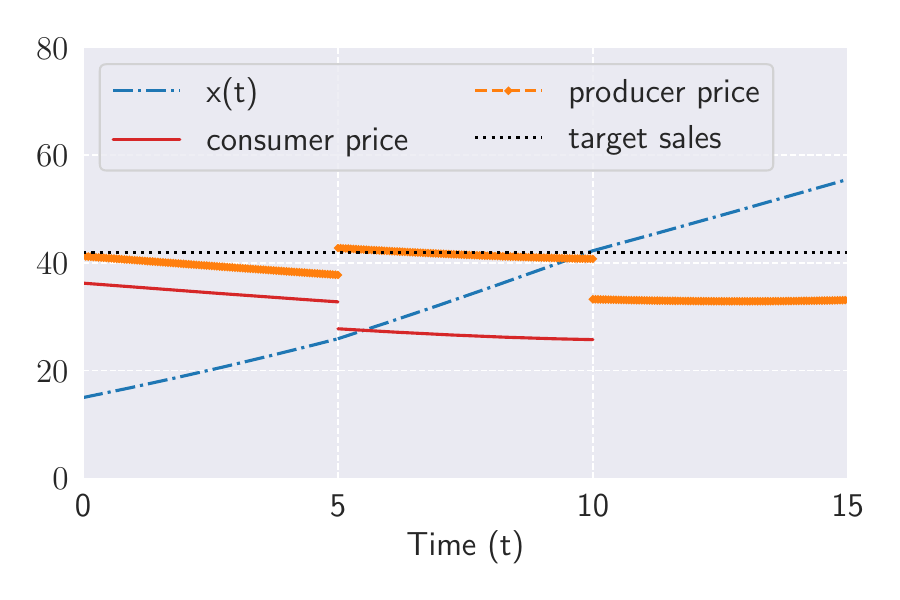}
            \caption{Target = 42}
            \label{fig:price_target_higher}
        \end{subfigure} \quad
\begin{subfigure}{0.35\linewidth}
            \centering
            \includegraphics[width=\textwidth]
          {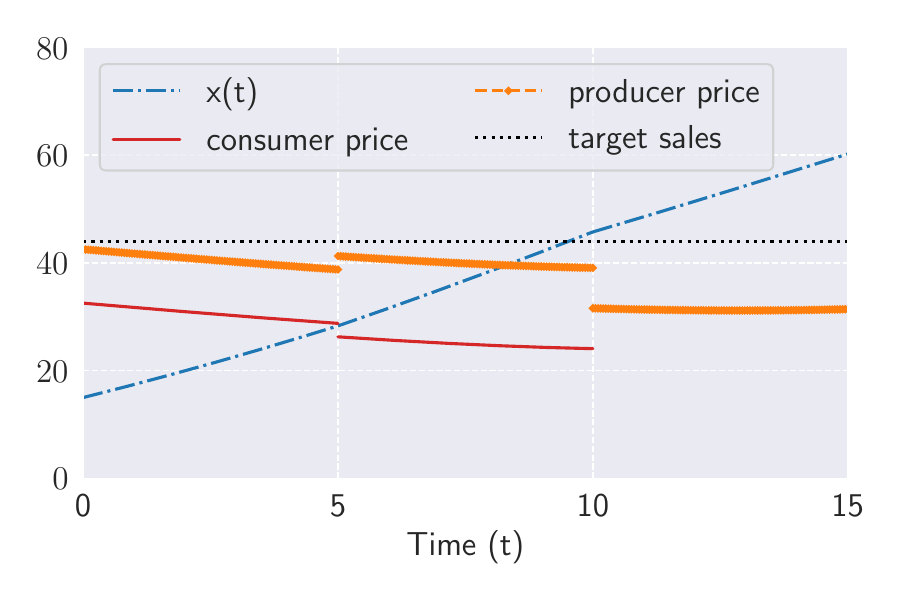}
            \caption{Target= 44}
            \label{fig:price_target_highest}
        \end{subfigure}
\caption{Cumulative sales, producer and consumer price when target is varied
from the benchmark case.}
\label{fig:prices_target}
\end{figure}
\begin{figure}[H]
\centering
\begin{subfigure}{0.22\linewidth}
            \centering
\includegraphics[width=\textwidth]{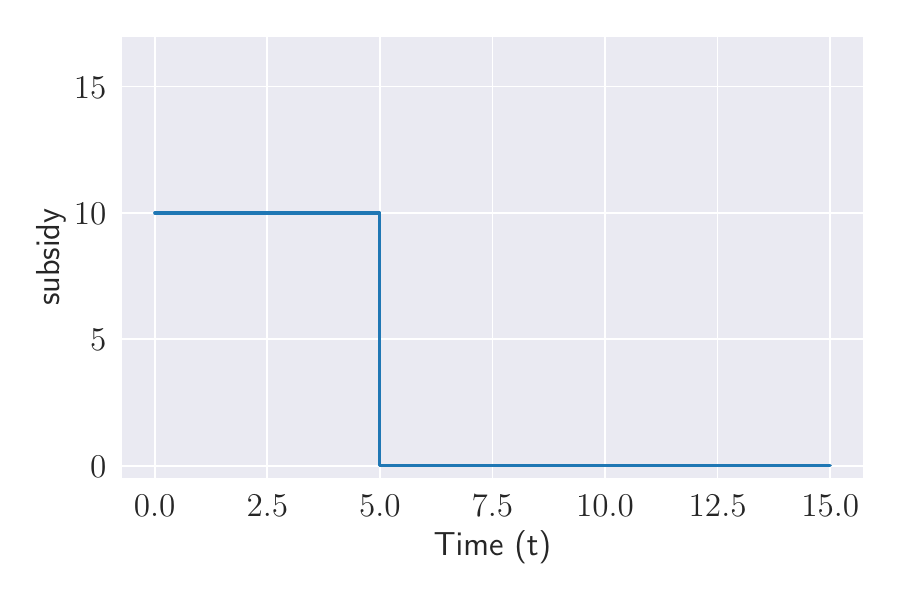}
            \caption{Target = 36}
            \label{fig:subsidy_target_lowest}
        \end{subfigure} 
\begin{subfigure}{0.22\linewidth}
            \centering
\includegraphics[width=\textwidth]{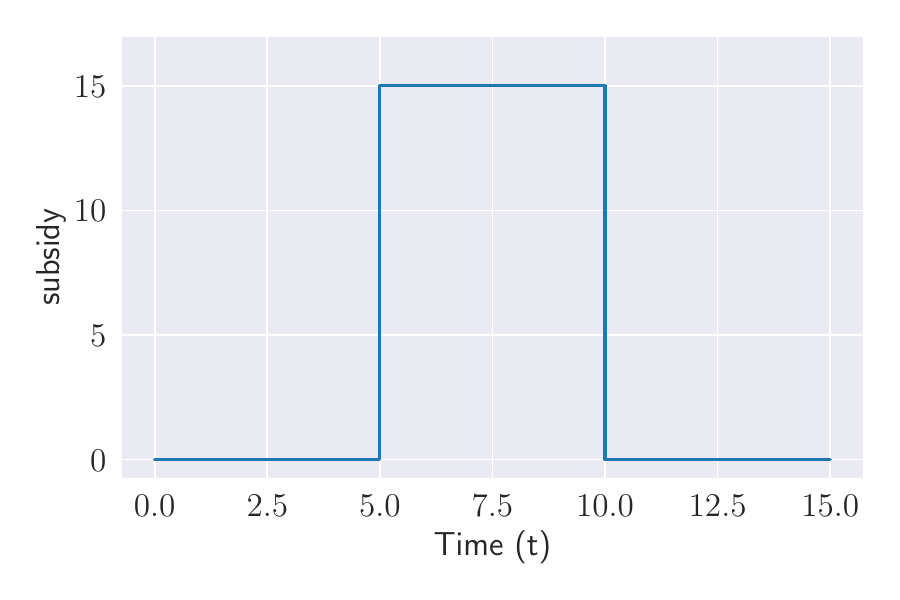}
            \caption{Target = 38}
            \label{fig:subsidy_target_lower}
        \end{subfigure} %
\begin{subfigure}{0.22\linewidth}
            \centering
\includegraphics[width=\textwidth]{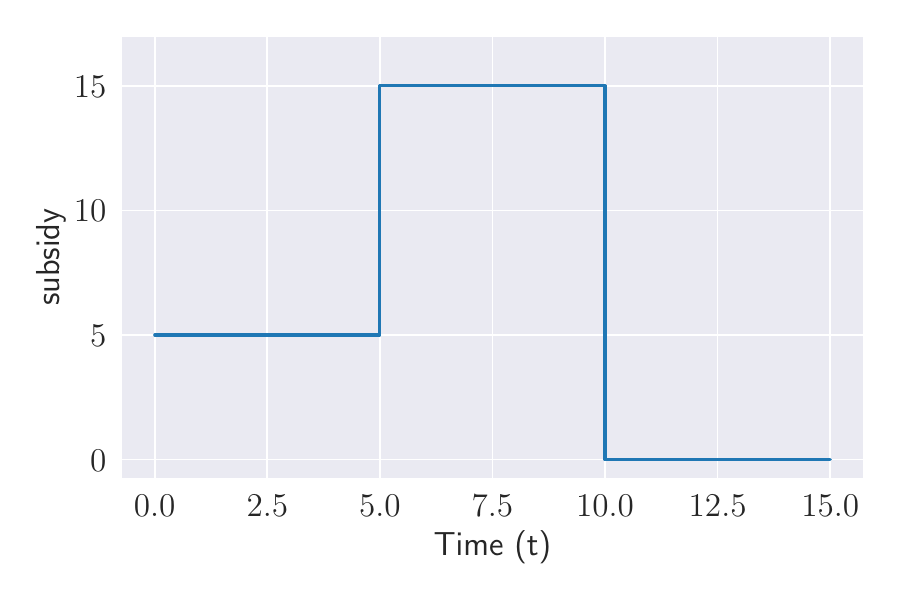}
            \caption{Target = 42}
            \label{fig:subsidy_target_higher}
        \end{subfigure} 
\begin{subfigure}{0.22\linewidth}
            \centering
            \includegraphics[width=\textwidth]
          {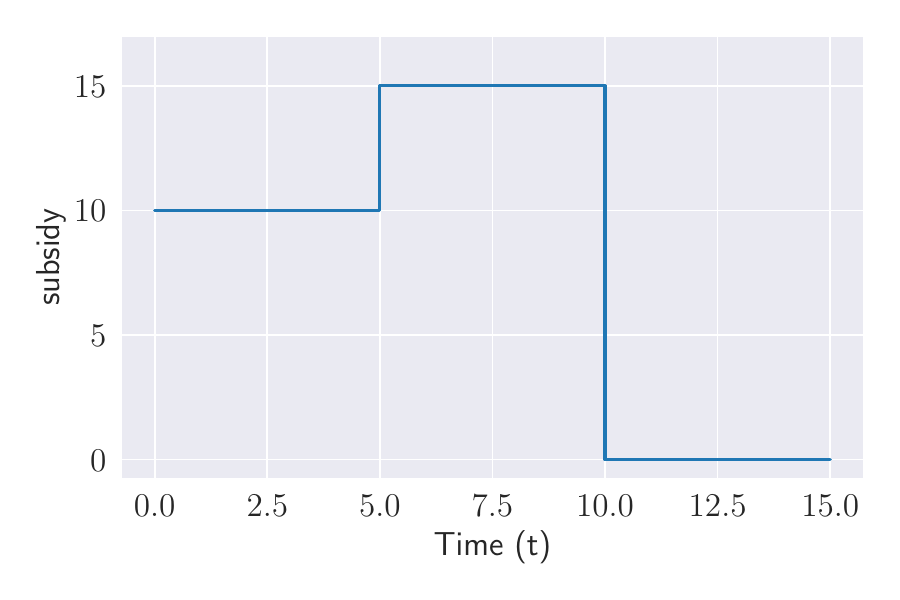}
            \caption{Target= 44}
            \label{fig:subsidy_target_highest}
        \end{subfigure}
\caption{{Equilibrium subsidy when target is varied from the
benchmark case.}}
\label{fig:subsidies_target}
\end{figure}

Consumers are paying a lower price in the subsidy scenarios and buying more.
Consequently, consumer surplus is higher with the subsidy. The firm is
clearly benefiting from the subsidy. Indeed, its profit is approximately $%
{970.88}$ with subsidy, and ${360.48}$ without a
subsidy, which is a ${170}\%$ increase. This big difference comes
from two sources. One is the higher sales volume, while the other is the
firm's strategic behavior discussed above. The environmental benefits are
more complicated to assess because they depend on a series of assumptions. A
first assessment can be done in terms of the number of gasoline cars
replaced by EVs, and the saving of gasoline consumption over the useful life
of an EV, which depends of the yearly driving distance by a car. From this
perspective, independent of how the computations are done, the conclusion
would be the same, i.e., subsidizing EVs reduces pollution emissions.
Ultimately, however, a comprehensive evaluation should consider all steps
involved in the production of the two types of cars from extraction of raw
materials to manufacturing and disposing of them when becoming obsolete.
Also, one should consider the sources of electricity used to feed the EVs.
Clearly, if the source is heavily polluting (coal, fuel), than the benefit
(if any) is much lower than when the electricity is produced with renewable
technologies.

\begin{figure}[H]
\centering
\begin{subfigure}{0.35\linewidth}
            \centering
\includegraphics[width=\textwidth]{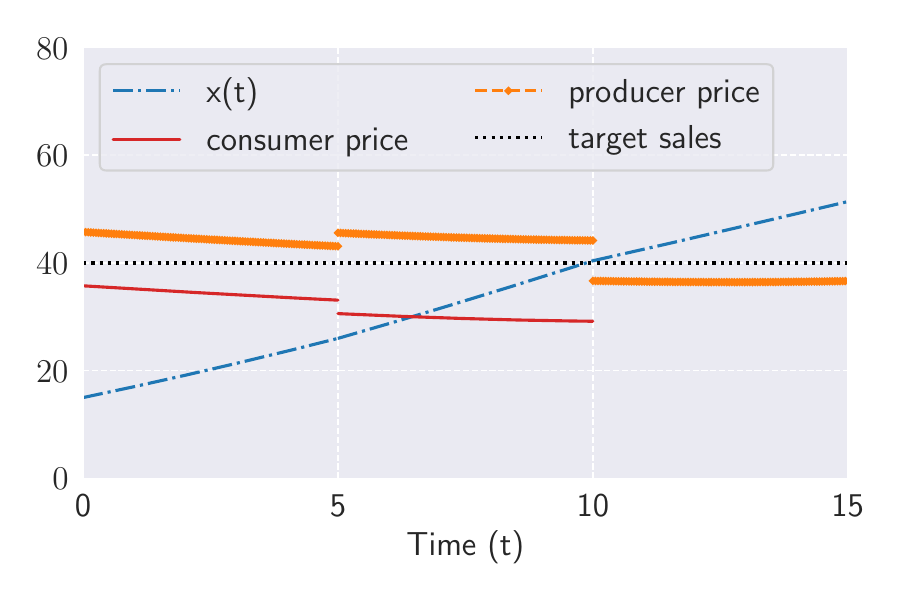}
            \caption{$b_2$= 0.72}
            \label{fig:_price_learn_lowest}
        \end{subfigure} 
\begin{subfigure}{0.35\linewidth}
            \centering
\includegraphics[width=\textwidth]{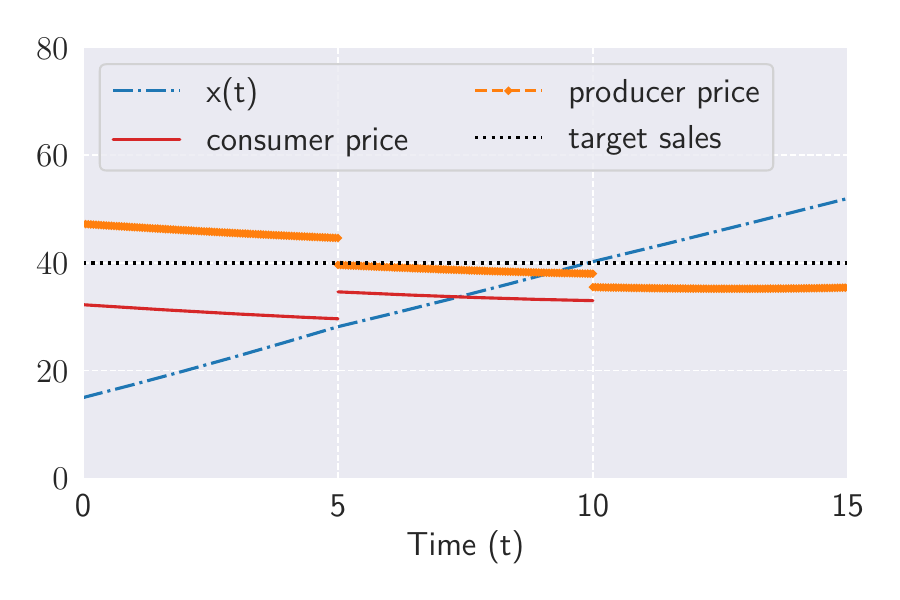}
            \caption{$b_2$= 0.76}
            \label{fig:price_learn_lower}
        \end{subfigure} \\
\begin{subfigure}{0.35\linewidth}
            \centering
\includegraphics[width=\textwidth]{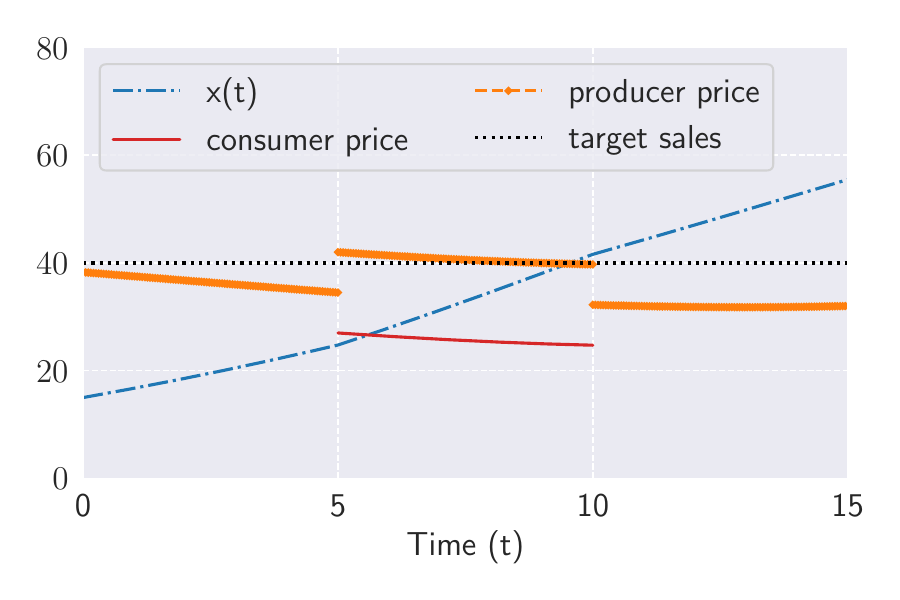}
            \caption{$b_2$= 0.84}
            \label{fig:price_learn_higher}
        \end{subfigure} 
\begin{subfigure}{0.35\linewidth}
            \centering
            \includegraphics[width=\textwidth]
          {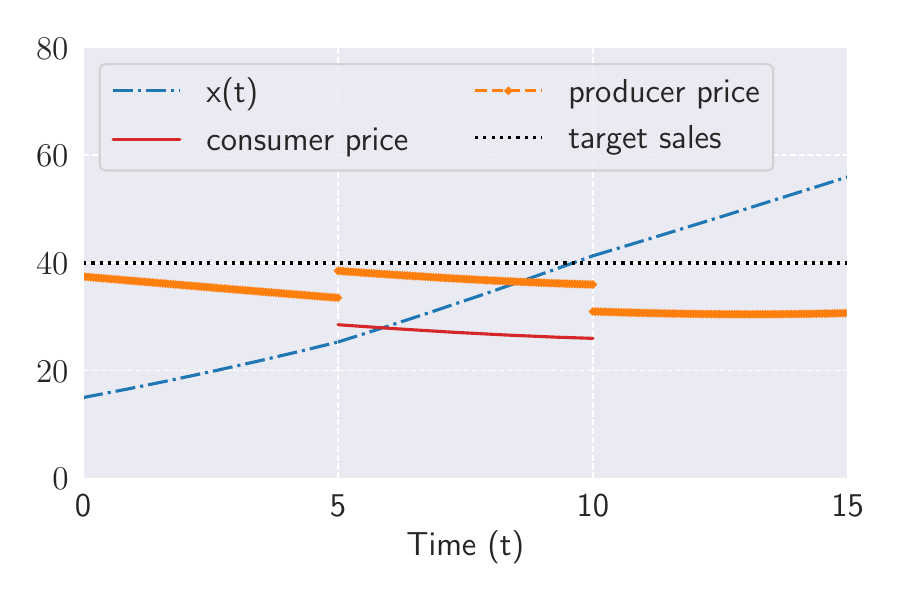}
            \caption{$b_2$= 0.88}
            \label{fig:price_learn_highest}
        \end{subfigure}
\caption{{Cumulative sales, producer and consumer price when $b_2$ is
varied from the benchmark case.}}
\label{fig:prices_learn}
\end{figure}

\subsection{Sensitivity analysis}

In this section, we vary the values of the main model's parameters and
assess the impact on the results.

\paragraph{Impact of the target value.}
We analyze two target values below and two above the benchmark level. As shown in Figure \ref{fig:prices_target}, reducing the target below the benchmark leads to a significant increase in consumer prices during the interval $[5,10]$. This increase corresponds to a reduction in the subsidy from $15$ in the benchmark case (Figure \ref{fig:subsidylevel}) to $0$ in the interval $[5,10]$ (Figure \ref{fig:subsidy_target_lowest}). For a target value of $38$, a lower subsidy is needed  in the interval $[0,5]$ compared to the benchmark case (Figure \ref{fig:subsidy_target_lower}). For target values above the benchmark, specifically $x_s = 42$ and $x_s = 44$, the subsidy is consistently equal to or higher than in the benchmark case across both time intervals $[0,5]$ and $[5,10]$, as illustrated in Figure \ref{fig:subsidy_target_higher}.
   The cost of the subsidy program for $%
x_{s}=36,38,42,44$ is ${800.77, 1090.93, 1593.29, 2127.73}$,
respectively. Note that the higher the target value, the higher the subsidy
and the government's cost, which is intuitive. What is less intuitive is the
order of magnitude in the changes. 
The cost of the subsidy program increases by ${165}$\%, when the
target is up by $22$\% (from $36$ to $44$). Finally, the higher the target,
the lower the after-subsidy price.

\begin{figure}[H]
\centering
\begin{subfigure}{0.22\linewidth}
            \centering
\includegraphics[width=\textwidth]{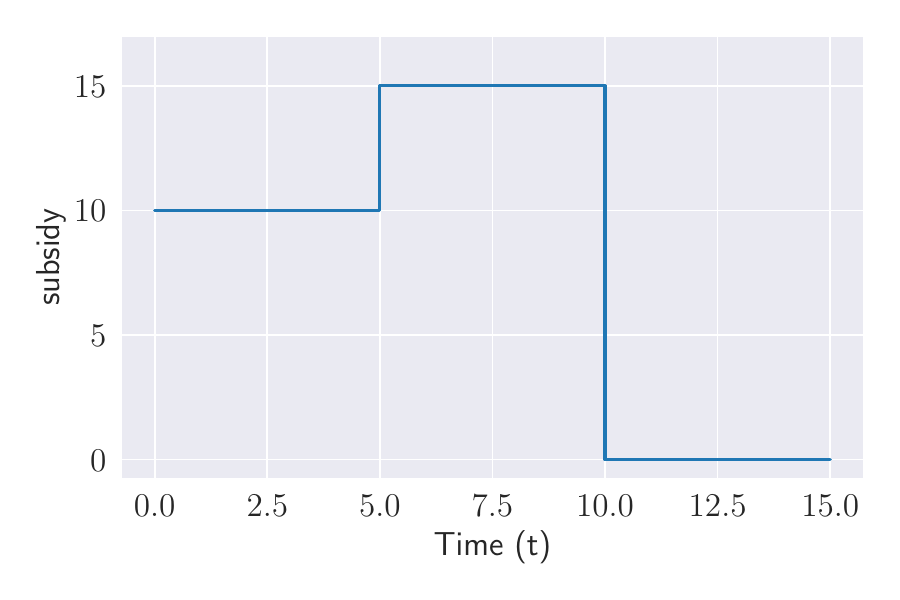}
            \caption{$b_2$= 0.72}
            \label{fig:subsidy_learn_lowest}
        \end{subfigure} 
\begin{subfigure}{0.22\linewidth}
            \centering
\includegraphics[width=\textwidth]{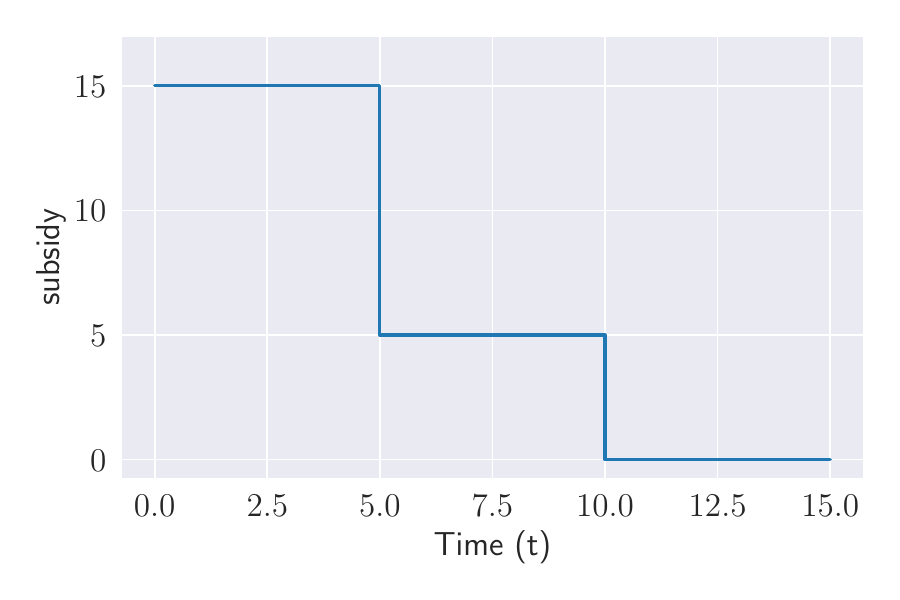}
            \caption{$b_2$= 0.76}
            \label{fig:subsidy_learn_lower}
        \end{subfigure} %
\begin{subfigure}{0.22\linewidth}
            \centering
\includegraphics[width=\textwidth]{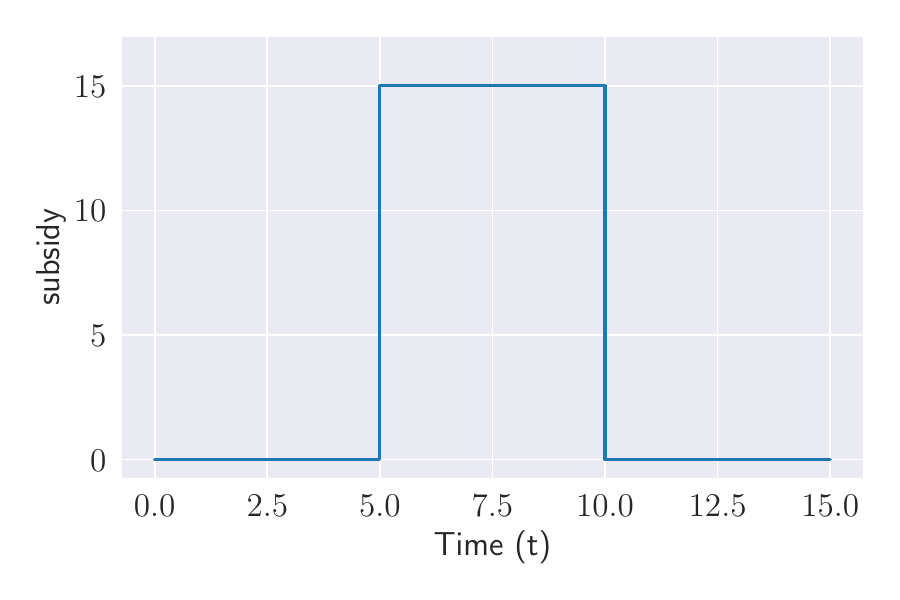}
            \caption{$b_2$= 0.84}
            \label{fig:subsidy_learn_higher}
        \end{subfigure} 
\begin{subfigure}{0.22\linewidth}
            \centering
            \includegraphics[width=\textwidth]
          {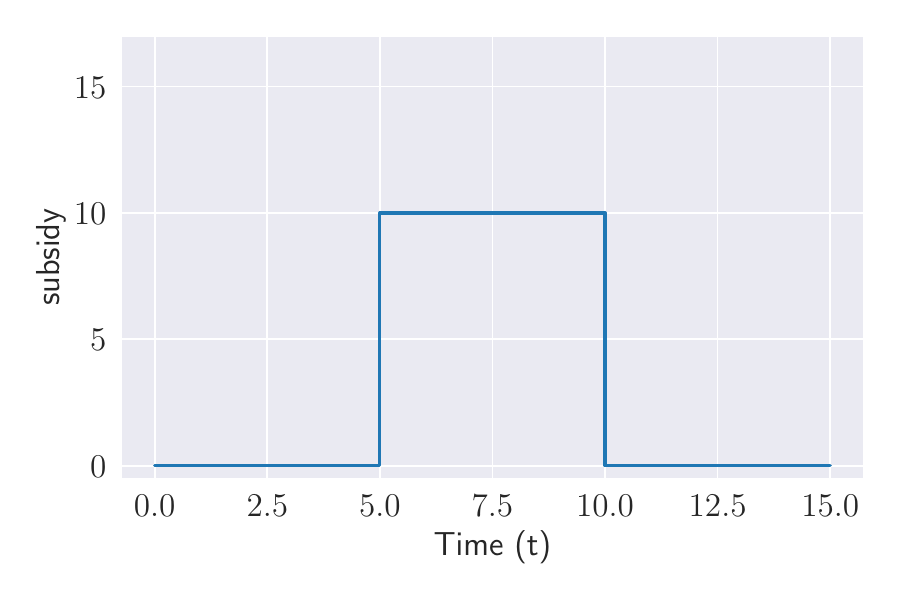}
            \caption{$b_2$= 0.88}
            \label{fig:subsidy_learn_highest}
        \end{subfigure}
\caption{{Equilibrium subsidy plan when $b_2$ is varied from the
benchmark case.}}
\label{fig:subsidies_learn}
\end{figure}
\begin{figure}[H]
\centering
\includegraphics[width=0.4\linewidth]{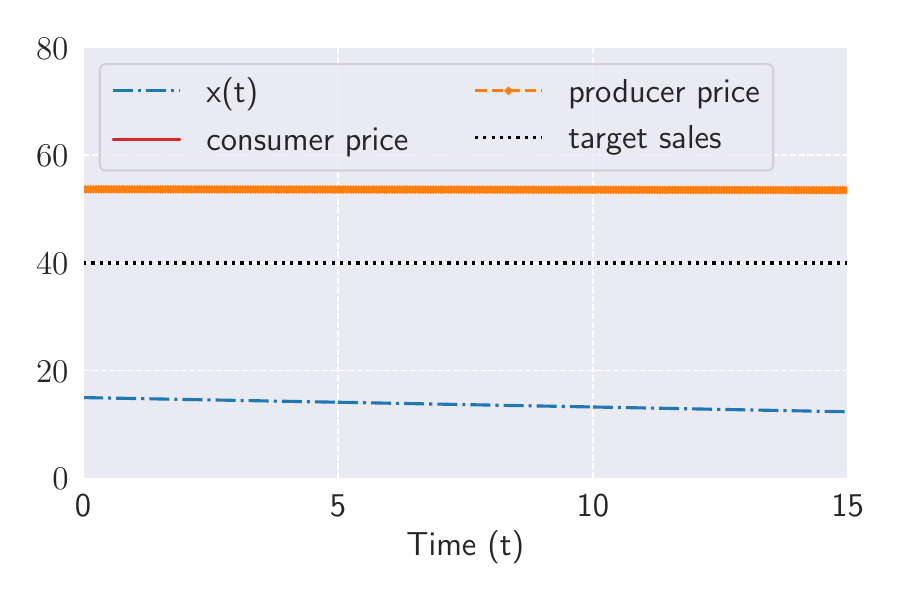}
\caption{{Cumulative sales, producer and consumer price when there is
no learning-by-doing ($b_2=0$).}}
\label{fig:price_no_lbd}
\end{figure}

\begin{figure}[H]
\centering
\begin{subfigure}{0.35\linewidth}
            \centering
\includegraphics[width=\textwidth]{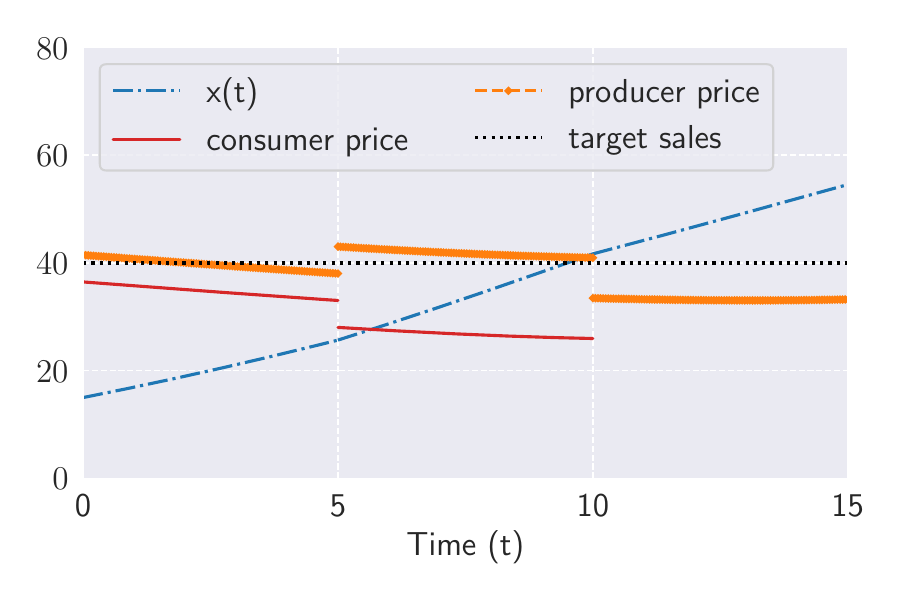}
            \caption{$\alpha_2$= 0.009}
            \label{fig:price_word_lowest}
        \end{subfigure} 
\begin{subfigure}{0.35\linewidth}
            \centering
\includegraphics[width=\textwidth]{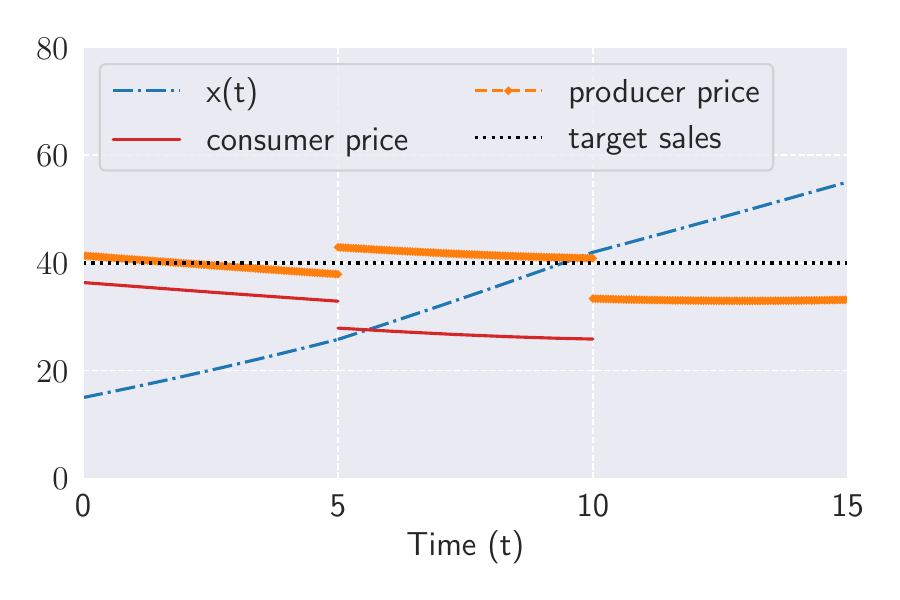}
            \caption{$\alpha_2$= 0.0095}
            \label{fig:price_word_lower}
        \end{subfigure} \\
\begin{subfigure}{0.35\linewidth}
            \centering
\includegraphics[width=\textwidth]{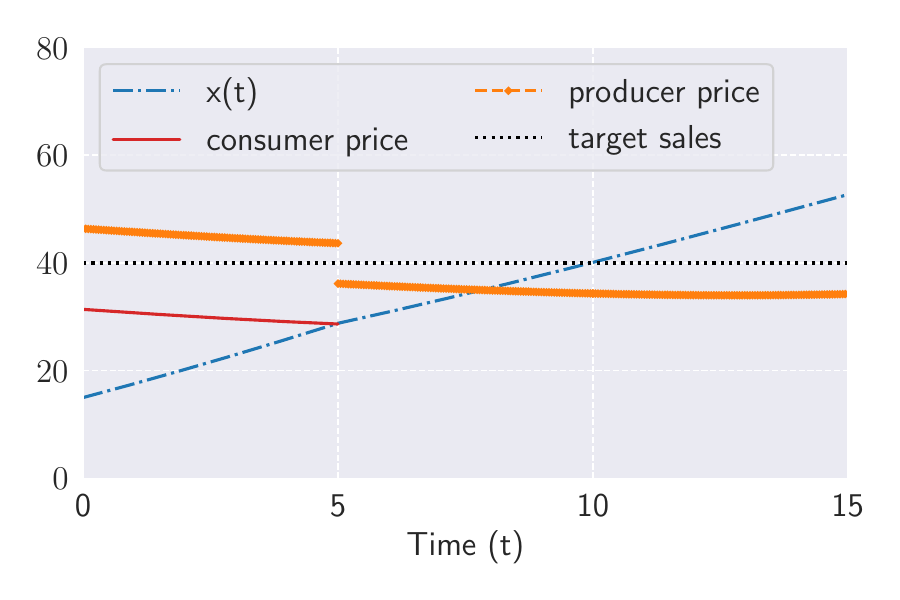}
            \caption{$\alpha_2$= 0.0105}
            \label{fig:price_word_higher}
        \end{subfigure} 
\begin{subfigure}{0.35\linewidth}
            \centering
            \includegraphics[width=\textwidth]
          {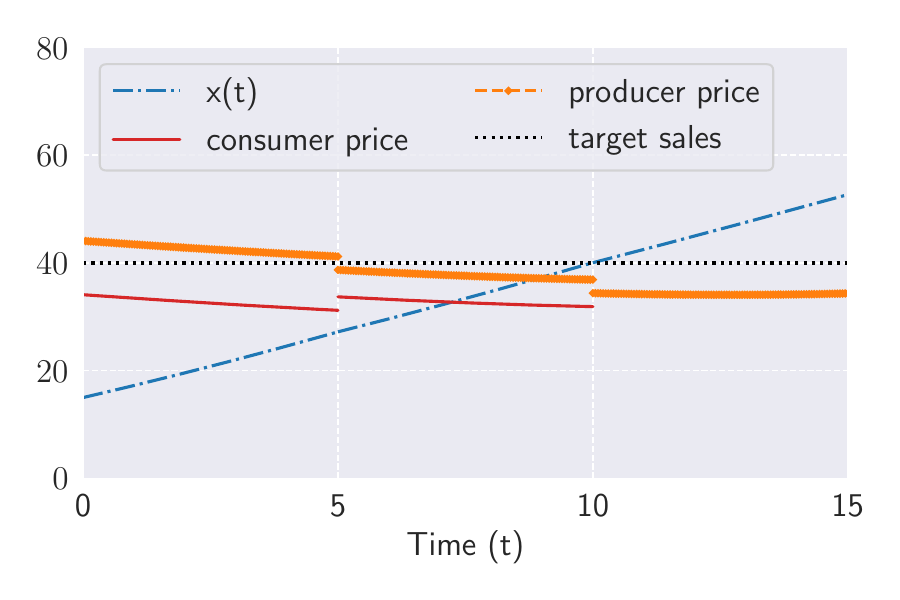}
            \caption{$\alpha_2$= 0.011}
            \label{fig:price_word_highest}
        \end{subfigure}
\caption{{Cumulative sales, producer and consumer price when $\protect%
\alpha_2$ is varied from the benchmark case.}}
\label{fig:prices_word}
\end{figure}
\paragraph{Impact of learning speed.}

As one can expect, a higher value of the learning speed $b_{2}$ only brings
good news. Indeed, we see in Figures \ref{fig:prices_learn} and \ref%
{fig:subsidies_learn} that increasing $b_{2}$ leads to lower price and
subsidy. Also, the cost to government is lower; for $%
b_{2}=0.72,0.76,0.84,0.88$, the subsidy cost is ${ 1957.39,
1644.83, 1159.79, 778.29}$, respectively. Clearly, the impact of the
learning speed is huge. Increasing $b_{2}$ by $22$\% (from 0.72 to 0.88)
cuts the subsidy budget by more than {$2.5$ times ($778.29$ when $%
b_{2}=0.88$ instead of $1957.39$ when $b_{2}=0.72$).}

{We can see in Figure \ref{fig:price_no_lbd} that if there is no
learning-by-doing, then government's subsidy program is not able to achieve
its target. This illustrates the complementary relationship between learning-by-doing and subsidy programs.}

\paragraph{Impact of word of mouth.}

\begin{figure}[H]
\centering
\begin{subfigure}{0.22\linewidth}
            \centering
\includegraphics[width=\textwidth]{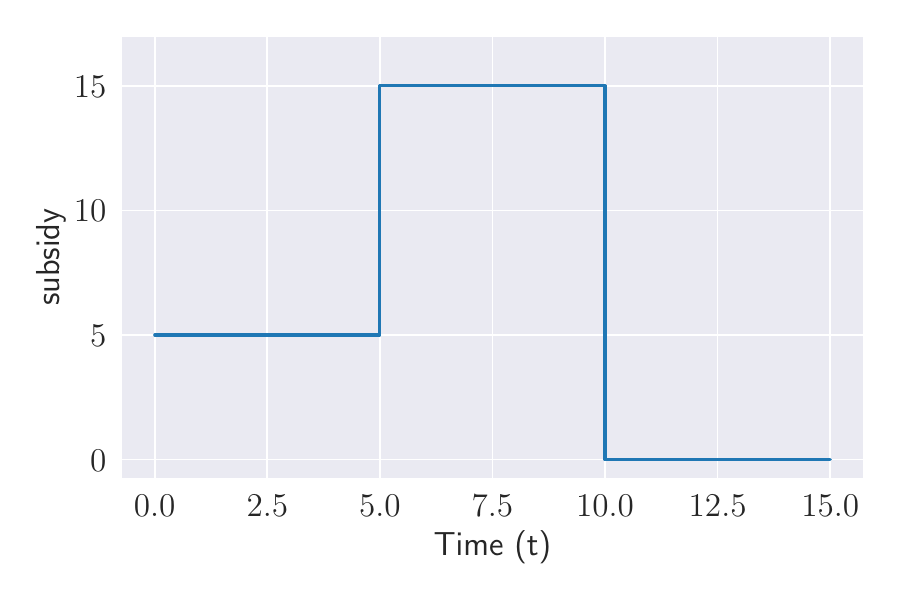}
            \caption{$\alpha_2$= 0.009}
            \label{fig:subsidy_word_lowest}
        \end{subfigure} 
\begin{subfigure}{0.22\linewidth}
            \centering
\includegraphics[width=\textwidth]{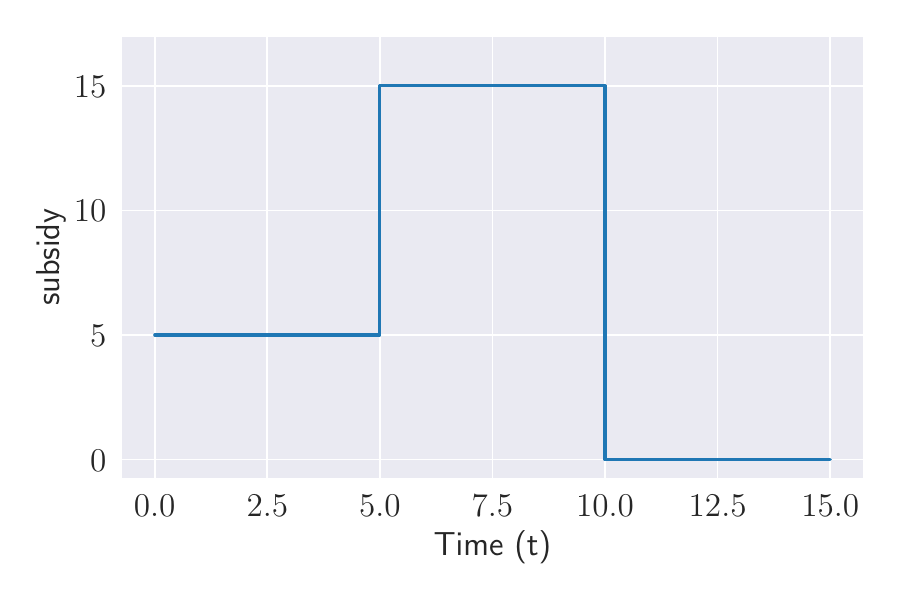}
            \caption{$\alpha_2$= 0.0095}
            \label{fig:subsidy_word_lower}
        \end{subfigure} %
\begin{subfigure}{0.22\linewidth}
            \centering
\includegraphics[width=\textwidth]{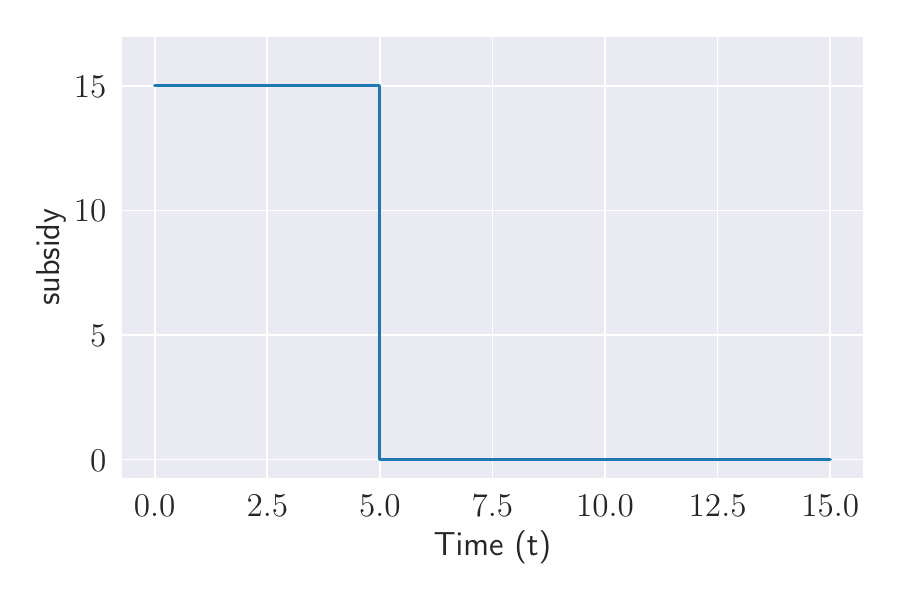}
            \caption{$\alpha_2$= 0.0105}
            \label{fig:subsidy_word_higher}
        \end{subfigure} 
\begin{subfigure}{0.22\linewidth}
            \centering
            \includegraphics[width=\textwidth]
          {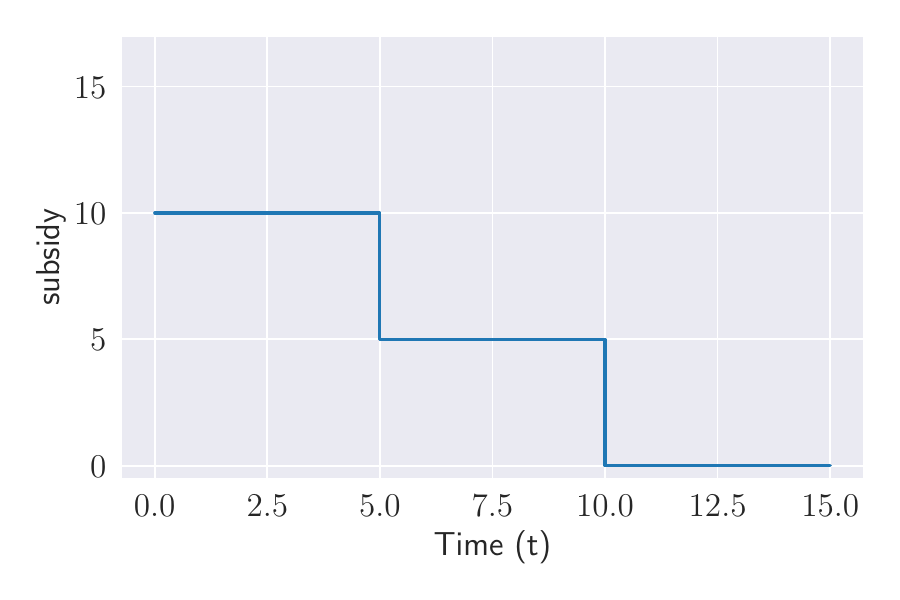}
            \caption{$\alpha_2$= 0.011}
            \label{fig:subsidy_word_highest}
        \end{subfigure}
\caption{{Equilibrium subsidy plan when $\protect\alpha_2$ is varied
from the benchmark case.}}
\label{fig:subsidies_word}
\end{figure}

We study the variation in equilibrium prices and subsidy with the change in $%
\alpha _{2}$ which measures the word-of-mouth effect. One explanation is that a higher $\alpha_2$
means larger market potential, which reduces the incentive to reduce the
price to boost demand. We can see in Figure %
\ref{fig:prices_word} as $\alpha _{2} $ increases, there is not much
variation in the price.  The subsidy required to reach the target decreases with increase in  $\alpha_2$ as illustrated in Figure %
\ref{fig:subsidies_word}.  The subsidy program cost for $\alpha
_{2}=0.009,0.0095,0.0105,0.011$ is given by ${%
1576.9, 1585.03, 1260.6, 1209.27}$, respectively. Here, increasing by 22\%
the value of $\alpha _{2}$ (from $0.009$ to $0.011$), leads to a 
{23\%} decrease in the budget.

\begin{figure}[H]
\centering
\begin{subfigure}{0.35\linewidth}
            \centering
\includegraphics[width=\textwidth]{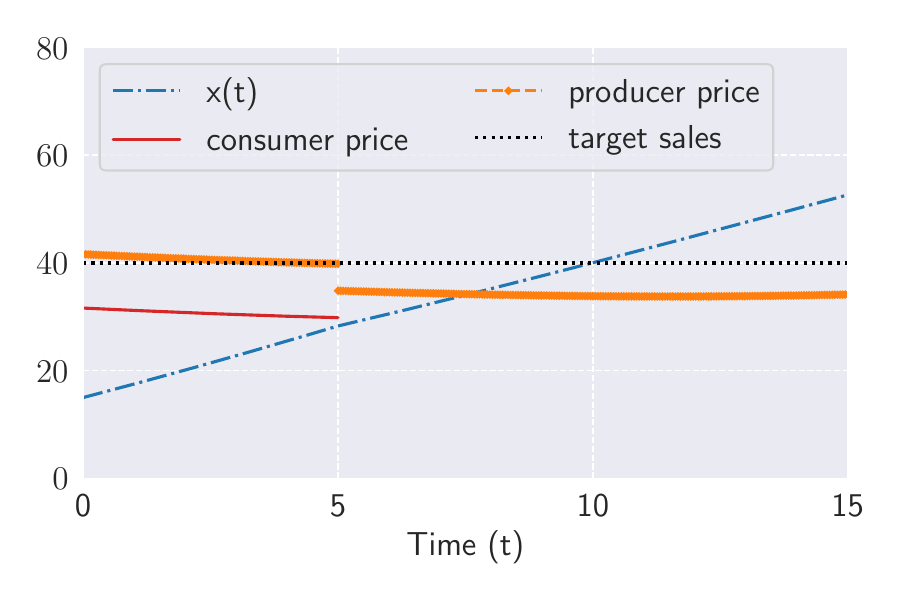}
            \caption{$\rho$= 0.06}
            \label{fig:rho6}
        \end{subfigure} 
\begin{subfigure}{0.35\linewidth}
            \centering
\includegraphics[width=\textwidth]{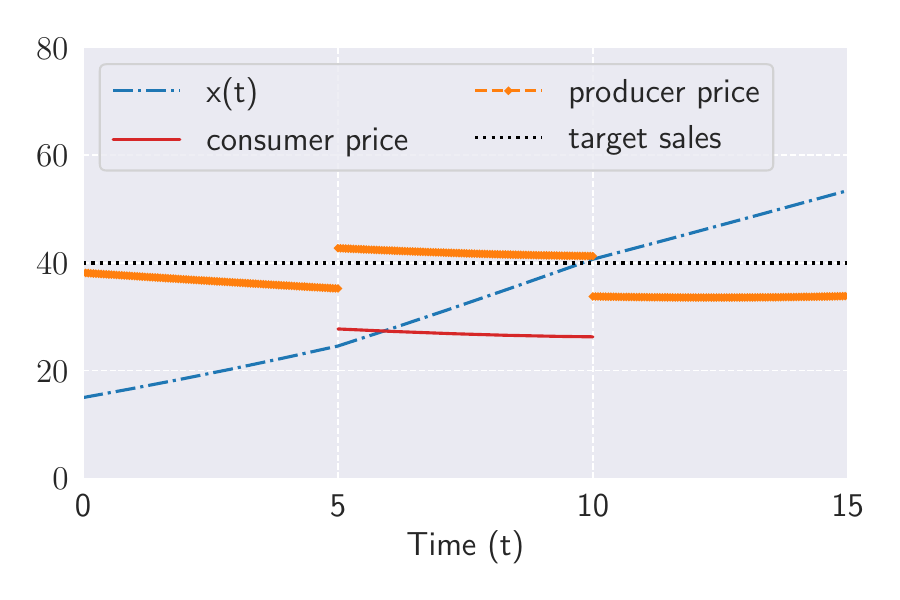}
            \caption{$\rho$= 0.08}
            \label{fig:rho8}
        \end{subfigure} \\
\begin{subfigure}{0.35\linewidth}
            \centering
\includegraphics[width=\textwidth]{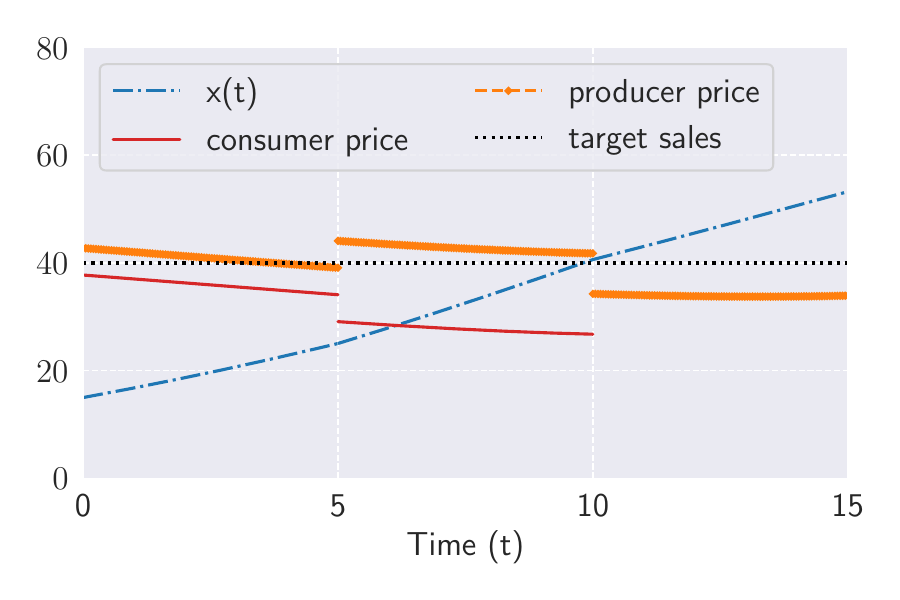}
            \caption{$\rho$= 0.12}
            \label{fig:rho12}
        \end{subfigure} 
\begin{subfigure}{0.35\linewidth}
            \centering
            \includegraphics[width=\textwidth]{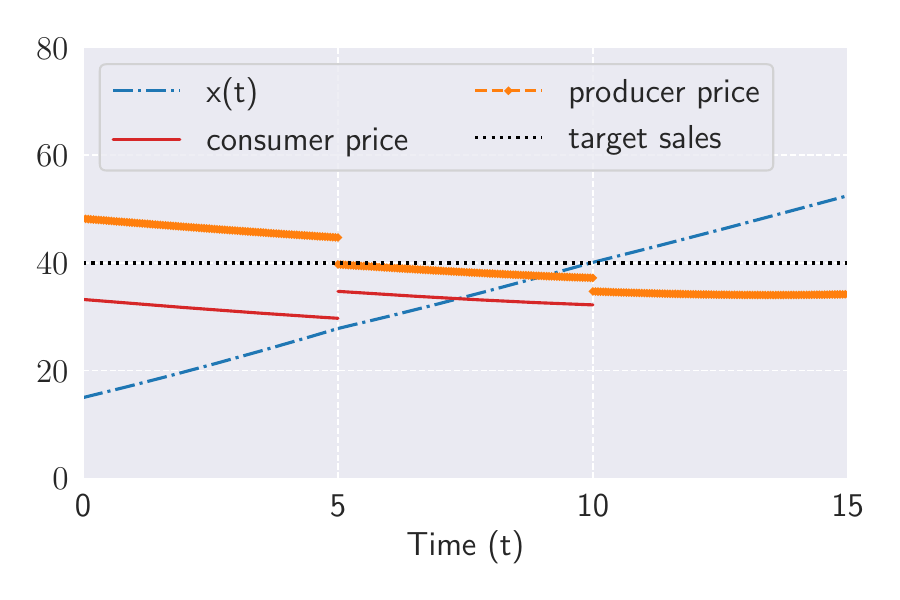}
            \caption{$\rho$= 0.14}
            \label{fig:rho14}
        \end{subfigure}
\caption{{Cumulative sales, producer and consumer price when $\protect%
\rho$ is varied from the benchmark case.}}
\label{fig:prices_rho}
\end{figure}

\begin{figure}[H]
\centering
\begin{subfigure}{0.24\linewidth}
            \centering
\includegraphics[width=\textwidth]{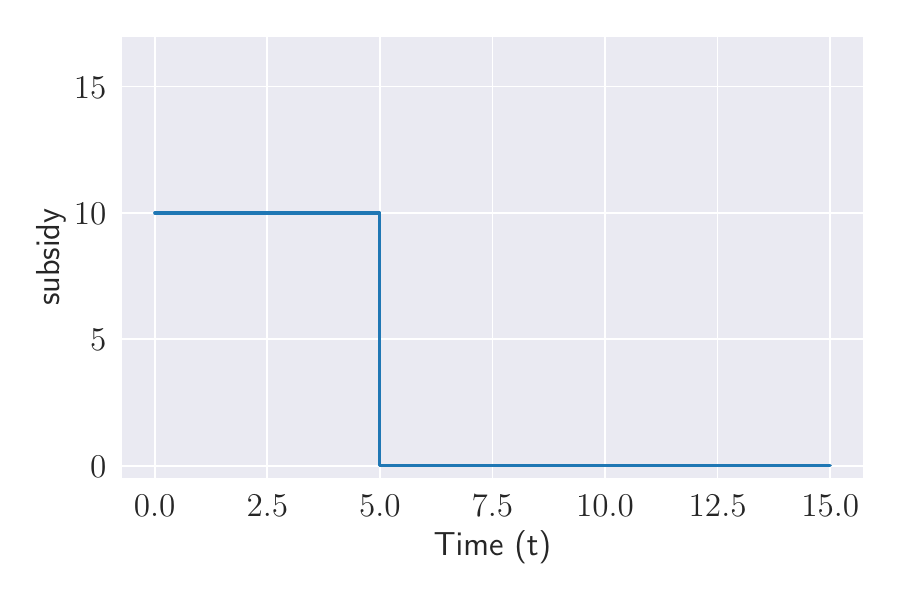}
            \caption{$\rho$= 0.06}
            \label{fig:rho6subs}
        \end{subfigure} 
\begin{subfigure}{0.24\linewidth}
            \centering
\includegraphics[width=\textwidth]{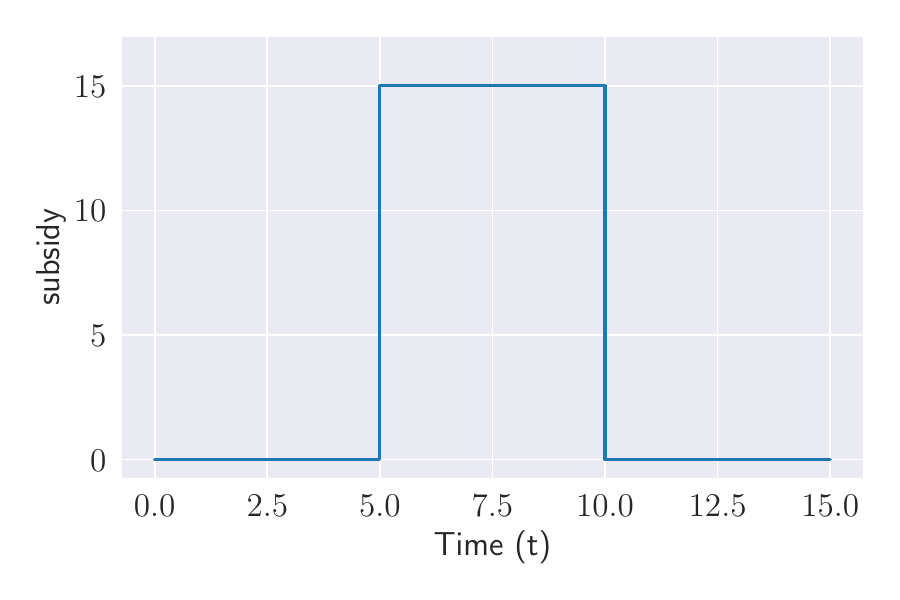}
            \caption{$\rho$= 0.08}
            \label{fig:rho8_subs}
        \end{subfigure} 
\begin{subfigure}{0.24\linewidth}
            \centering
\includegraphics[width=\textwidth]{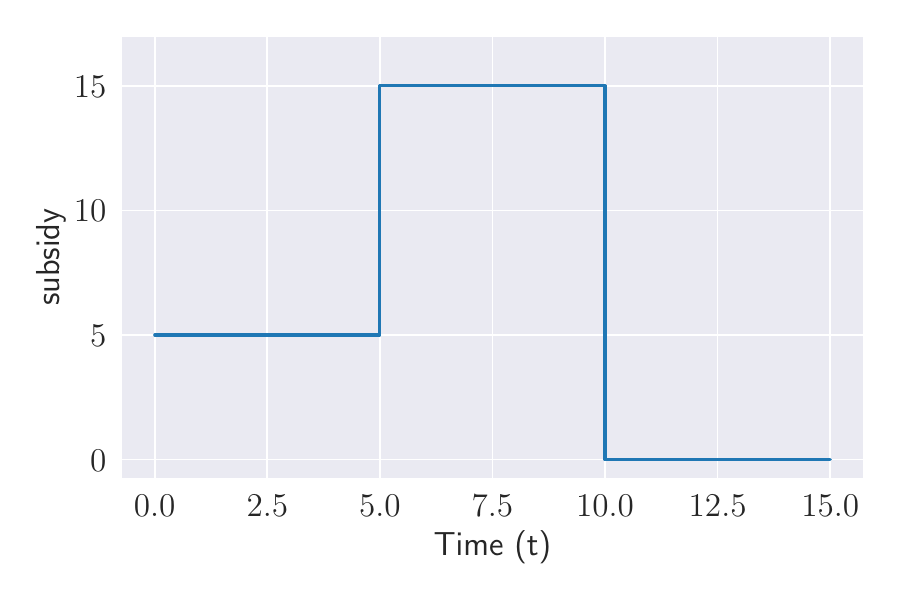}
            \caption{$\rho$= 0.12}
            \label{fig:rho12subs}
        \end{subfigure} 
\begin{subfigure}{0.24\linewidth}
            \centering
            \includegraphics[width=\textwidth]{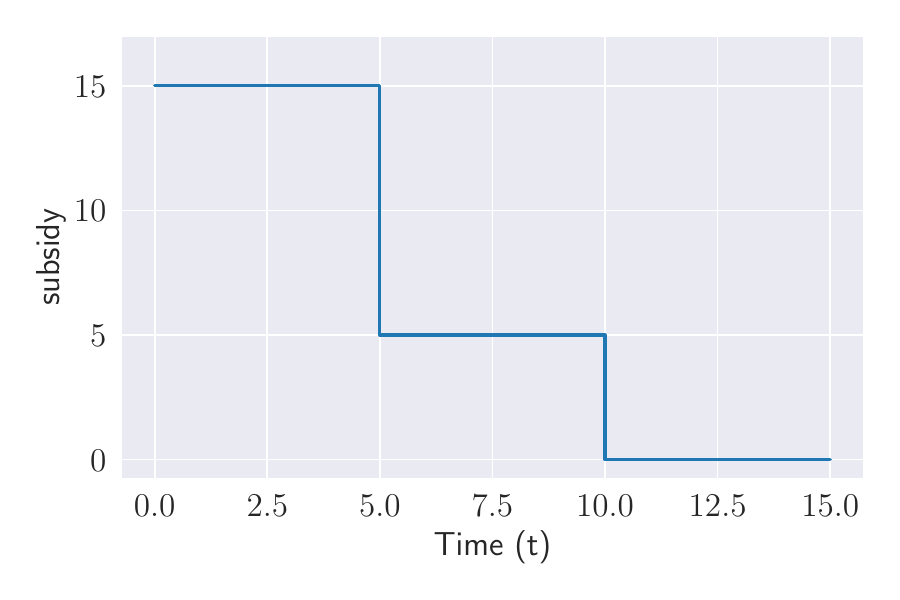}
            \caption{$\rho$= 0.14}
            \label{fig:rho14_subs}
        \end{subfigure}
\caption{{Equilibrium subsidy when $\rho$ is varied from the
benchmark case.}}
\label{fig:subsidy_rho}
\end{figure}

\paragraph{Impact of the discount factor} {In this scenario, we expand the variation in the discount factor up to
40\%, compared to the previously considered 10\% for other parameters, to
observe more significant impacts on both producer and consumer prices, see Figure \ref{fig:prices_rho}.
Relative to the benchmark case, a lower discount factor leads to a decline
in the producer price, as well as a reduction in the government's subsidy,
since less immediate intervention is needed to achieve its target. In
contrast, when the future is discounted more aggressively than in the
benchmark case, the government responds by increasing the subsidy to
maintain its policy goals, which subsequently drives the producer price even
higher, see Figure \ref{fig:subsidy_rho}. The cost of the subsidy program for $\rho =0.06,0.08,0.12,0.14$ is $%
929.89,1328.57,1362.38,1410.43$, respectively. }

\paragraph{Impact of subsidy adjustments.}

\begin{figure}[H]
\centering
\begin{subfigure}{0.45\linewidth}
            \centering
\includegraphics[width=\textwidth]{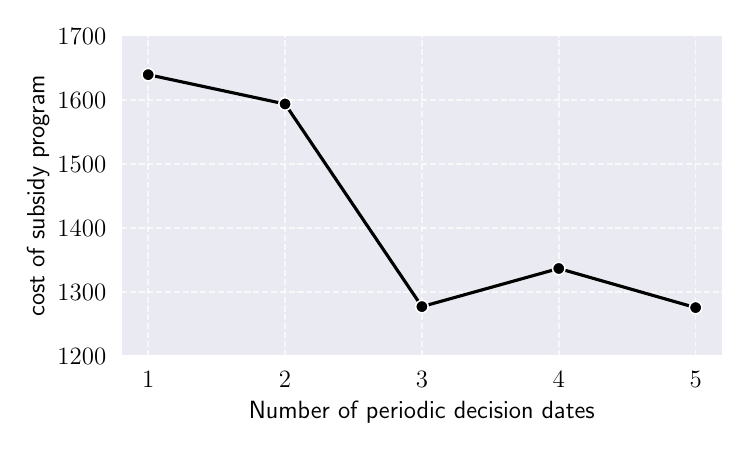}
            \caption{$C=10$}
            \label{fig:subsidy_cost_50}
        \end{subfigure} 
\begin{subfigure}{0.45\linewidth}
            \centering
\includegraphics[width=\textwidth]{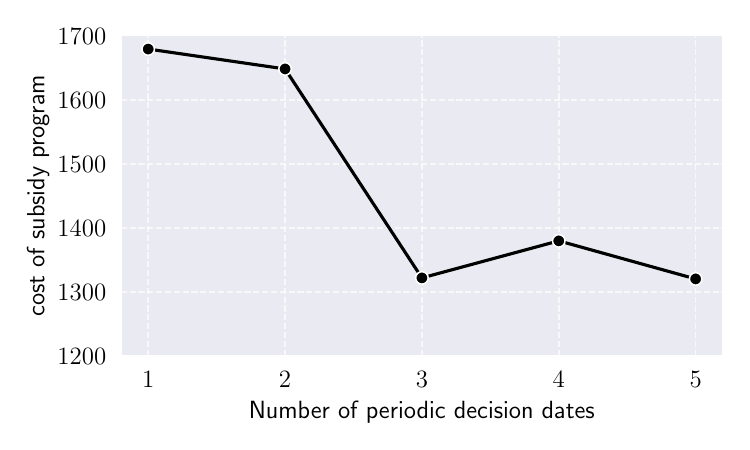}
            \caption{$C=50$}
            \label{fig:subsidy_cost_100}
        \end{subfigure}
\caption{{Variation in the cost of the subsidy program with number of
periodic decision dates}}
\label{fig:subsidies_cost}
\end{figure}

\begin{figure}[t]
\centering
\begin{subfigure}{0.45\linewidth}
            \centering
\includegraphics[width=\textwidth]{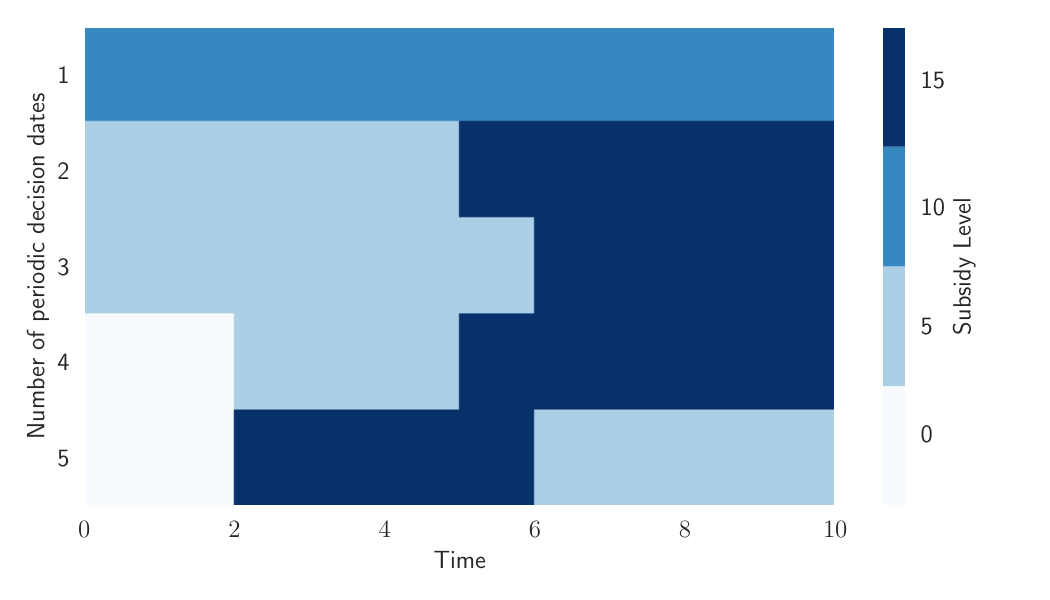}
            \caption{$C=10$}
            \label{fig:subsidy_time_50}
        \end{subfigure} 
\begin{subfigure}{0.45\linewidth}
            \centering
\includegraphics[width=\textwidth]{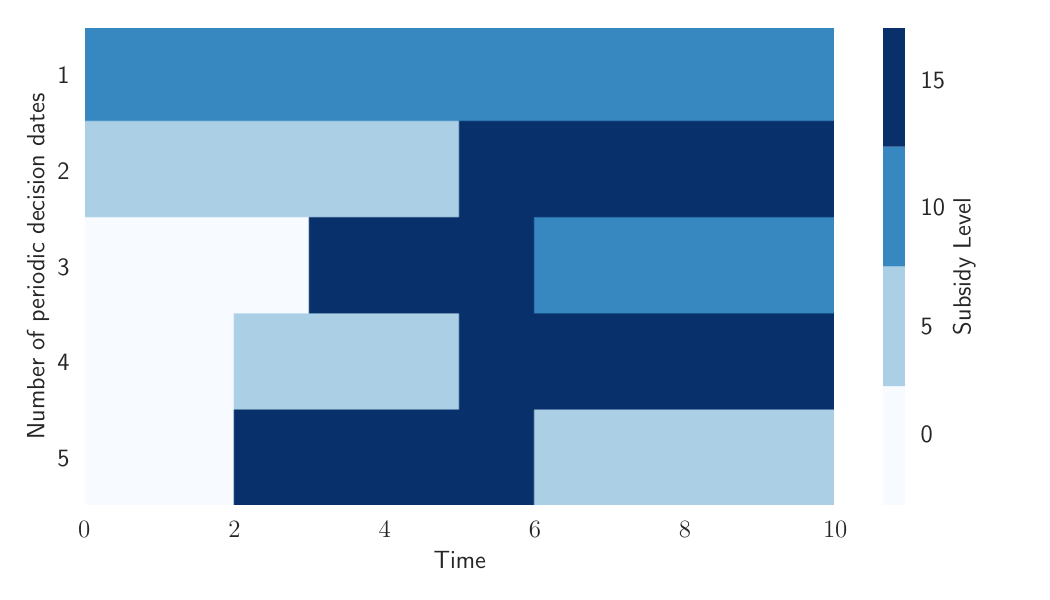}
            \caption{$C=50$}
            \label{fig:subsidy_time_100}
        \end{subfigure}
\caption{Variation in subsidy levels over time for different number of
periodic decision dates.}
\label{fig:subsidies_time}
\end{figure}

{Finally, we vary the number of subsidy adjustments in the time
interval $[0,10]$. As the number of adjustments increase in {Figure} \ref%
{fig:prices_dec}, the terminal price remains almost the same.  In Figure \ref{fig:subsidies_cost},
we plot the variation in the cost of the subsidy program with the number of
periodic decision dates. The non-monotonicity is due to the fact that the
firm can adjust the price to the subsidy policy of the government, which in
turn affects the cost of the subsidy program.} A larger number of changes
give more degrees of freedom in adjusting the subsidy levels to reach the
target. Therefore, the ultimate impact will depend on the target and the
fixed cost of each change in the subsidy. {When the fixed cost of
modifying the subsidy increases from $10$ (Figure \ref{fig:subsidy_time_50})
to $50$ (Figure \ref{fig:subsidy_time_100}), the government may opt to delay
offering the subsidy until electric vehicle sales begin to grow--at which
point higher production volumes drive down costs and reduce the total
subsidy cost. }

\section{Conclusions}

\label{sec:conclusion} In this paper, we provide a verification theorem to
characterize the feedback-Stackelberg equilibrium in a differential game
between a government and a firm. While the firm acts at each instant of
time, the government intervenes only at certain discrete time instants to
adjust the subsidy level. To the best of our knowledge, it is the \USmodified{first} time
that a feedback-Stackelberg equilibrium is determined in a differential game
with one player using impulse control. Also, it is the first paper in the
diffusion models literature that implements discrete changes to the subsidy,
which is more realistic than assuming a continuous modification of its
level. 

\begin{figure}[t]
\centering
\begin{subfigure}{0.35\linewidth}
            \centering
\includegraphics[width=\textwidth]{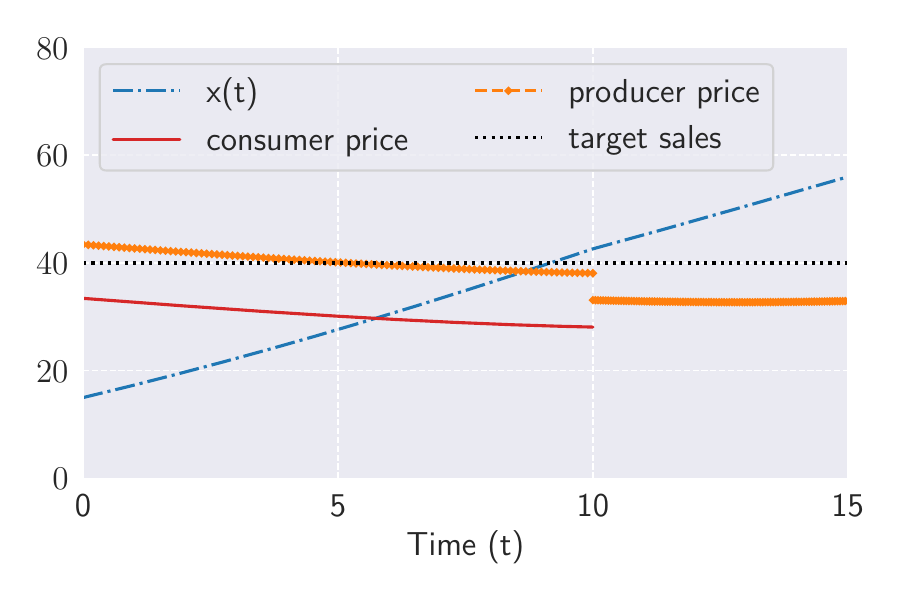}
            \caption{$N$= 1}
            \label{fig:price_dec_1}
        \end{subfigure} 
\begin{subfigure}{0.35\linewidth}
            \centering
\includegraphics[width=\textwidth]{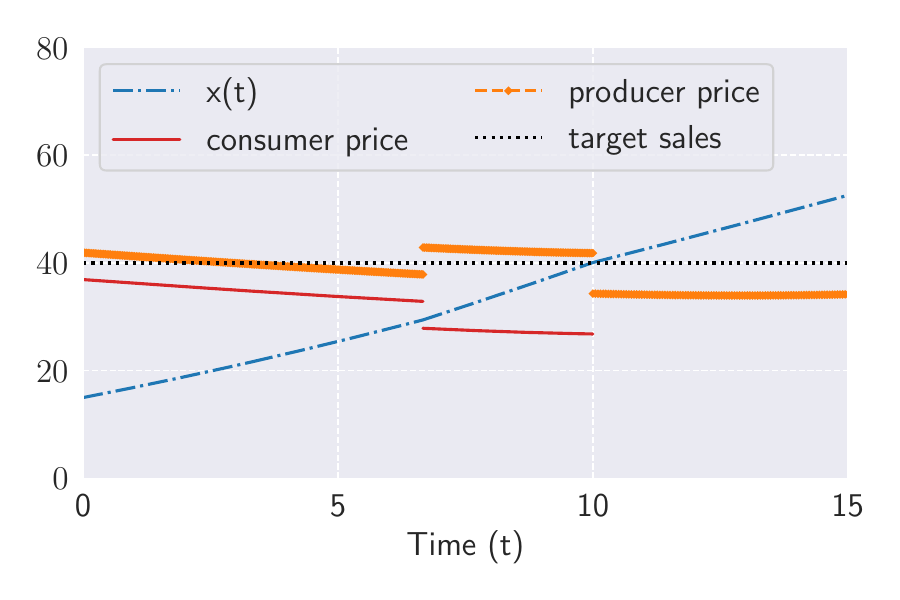}
            \caption{$N=3$}
            \label{fig:price_dec_3}
        \end{subfigure} \\
\begin{subfigure}{0.35\linewidth}
            \centering
\includegraphics[width=\textwidth]{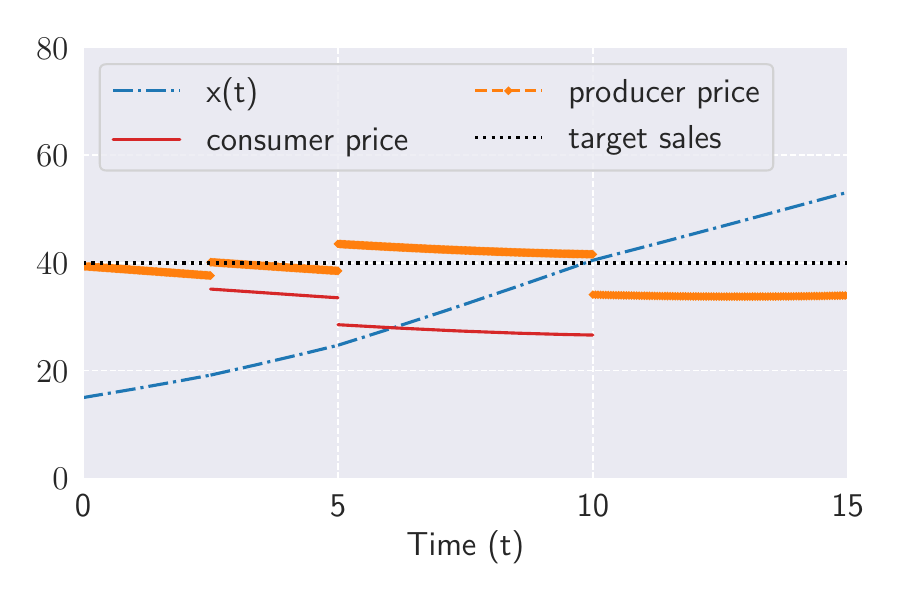}
            \caption{$N=4$}
            \label{fig:fig:price_dec_4}
        \end{subfigure} 
\begin{subfigure}{0.35\linewidth}
            \centering
            \includegraphics[width=\textwidth]
          {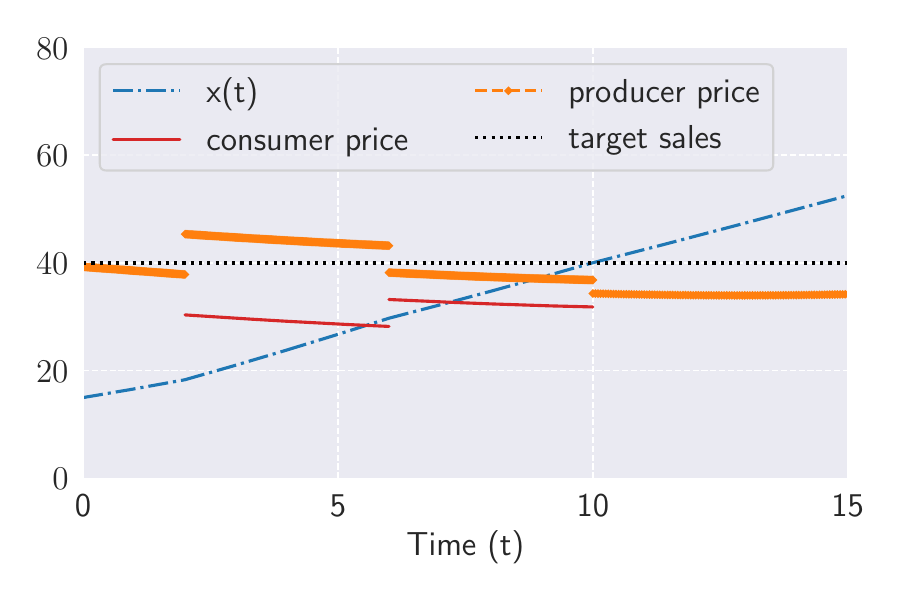}
            \caption{$N=5$}
            \label{fig:fig:price_dec_5}
        \end{subfigure}
\caption{{Cumulative sales, producer and consumer price when $\protect%
\rho$ is varied from the benchmark case.}}
\label{fig:prices_dec}
\end{figure}

\subsection{\USmodified{Managerial Contributions}}

In summary, our contribution is twofold. 

\begin{enumerate}
\item We provide a tool to decision makers in public sector to design a
subsidy program that accounts for the firm's behavior, cost dynamics
(learning-by-doing effect), and cumulative-sales target. 

\item We give some interesting insights into the equilibrium strategies and
how they vary with the parameter values. The main takeaways here are: 

\begin{enumerate}
\item Consumers benefit from a subsidy program, that is, they pay less than
what they would have paid in the absence of such program. 

\item The monopolist takes advantage of the subsidy program by raising its
price. How such behavior could be deterred is not an easy question in a
context where governmental price control is not on the menu. (More on this
below.) 

\item The presence of learning-by-doing seems to be a must for establishing
a subsidy program. Our results seem to show than in the absence of cost
decline, a subsidy program is useless, and the government should not engage
in such a program. 

\item The fixed cost incurred by the government each time it changes the
subsidy determines the frequency of such changes. Although such result is
expected, our model allows to experiment with different frequencies and
assess their impact on the realization of the target and on the cost of the
program.
\end{enumerate}
\end{enumerate}

\subsection{\USmodified{Extensions}}

In any analytical work, one has to make some delicate trade-offs between the
realism of the model, without including all details, and solvability. Our
work does not escape  this rule and we believe it is of intellectual and
practical interests to consider the following extensions:

\begin{enumerate}
\item We assumed that the cost function is affine in cumulative production,
which may be a good approximation in the short run, but not in the long
term. Based on empirical work, \cite{Levitt2013} concludes that
\textquotedblleft . . . learning is nonlinear: large gains are realized
quickly, but the speed of progress slows over time.\textquotedblright\ \ One
option is to adopt a hyperbolic cost function that allows to capture this
nonlinearity and insure that the cost remains positive for any level of
cumulative production; see, e.g., \cite{Janssens2014}. However,
such modification comes with the additional difficulty in determining the
feedback-Stackelberg equilibrium as the game would not be anymore
linear-quadratic {and one would need to numerically solve the HJB partial
differential equation}. \ 

\item Another extension is to let the timing of subsidy adjustments be also
chosen optimally. This extension would require challenging methodological
developments as we do not dispose yet of a theorem characterizing the
equilibrium in such setup.

\item As mentioned above, the monopolist sets its price at a higher value in
the presence of a subsidy program. An interesting question is to verify if
this behavior still holds when we increase the number of firms. The analysis
of a duopoly could already give some indications with this respect.

\item An attempt to estimate empirically (at least) some of the parameter
values would be clearly of interest.
\end{enumerate}

\appendix
\renewcommand{\thesection}{\Alph{section}}
\renewcommand{\thesubsection}{\Alph{section}.\arabic{subsection}}
\USmodified{\section{Proofs}}

\subsection{Proof for the analytical expression of $k_2(t)$}

\label{appen_k2} $k_2(t)$ satisfies the Riccati differential equation from $t\in [0, T]$ given below:
\begin{equation}
\rho k_2(t) - \dot{k}_2(t) = \frac{\beta}{2}\left(\frac{w_2}{\beta} +
k_2(t)\right)^2.  \label{eq:k2app}
\end{equation}
Define $u(t) = w_2/\beta + k_2(t)$ so that $k_2(t) = u(t)- w_2/\beta$.
Then \eqref{eq:k2app} becomes $
\rho\Bigl(u(t) - w_2/\beta\Bigr) - \dot{u}(t) = (\beta/2) u(t)^2$
or equivalently, we can write the standard Riccati differential equation with constant coefficients as follows:
\begin{align}
    \dot{u}(t) = \rho\, u(t) - \rho w_2/\beta - (\beta/2)\,
u(t)^2. \label{eq:riccati_st}
\end{align}

Case 1: Suppose that $\rho^2 - 2\rho\,w_2\geq 0$ and let $\Delta := \sqrt{\rho^2
- 2\rho\,w_2}$, then we obtain the Riccati equation: $\dot{u}(t) = -(\beta/2) \left((u(t) - \rho/\beta)^2 -
(\Delta/\beta)^2 \right).$
On rearranging, we obtain: 
\begin{equation*}
\frac{d u(t)}{(u(t) - \rho/\beta)^2-(\Delta/\beta)^2} = -\frac{\beta}{2} dt.
\end{equation*}
Next, we substitute $m(t)=\frac{u(t)-\rho/\beta}{\Delta/\beta}$ in the above equation. Then, we
obtain: $d m(t)/(1-m(t)^2) = (\Delta/2) dt.$
Taking integral on both sides, we obtain: $
\tanh^{-1} \left( m(t)\right) = (\Delta/2)(t+\zeta_\Delta).$
We substitute $m(t)=\frac{u(t)-\rho/\beta}{\Delta/\beta}$ to obtain 
$u(t)=(\rho/\beta) + (\Delta/\beta)\tanh\left((\Delta/2)(t+\zeta_\Delta)\right).$
On substituting $u(t)=w_2/\beta +k_2(t)$, we get: 
\begin{equation*}
k_2(t)= \frac{\rho -w_2+ \Delta\,\tanh\left((\Delta/2)%
(t+\zeta_\Delta)\right)}{\beta}.
\end{equation*}
$\zeta_\Delta$ is determined by the terminal condition, $k_2(T)=0$, so that \[(\rho-w_2)/\beta + (\Delta/\beta)\,\tanh\!\left((\Delta/2)(T+\zeta_\Delta)\right)=0.\] 
Multiplying both sides by $\beta$ and
solving for the hyperbolic tangent term, we obtain $
\tanh\!\left((\Delta/2)(T+\zeta_\Delta)\right) = -(\rho - w_2)/
\Delta$.
Taking the inverse hyperbolic tangent (i.e., $\operatorname{arctanh}$) on both
sides, $(\Delta/2)(T+\zeta_\Delta) = \operatorname{arctanh}\!\left((w_2-\rho)/\Delta\right).$
Thus, the constant $\zeta_\Delta$ is given by:
\begin{equation*}
\zeta_\Delta= \frac{2}{\Delta}\operatorname{arctanh}\!\left(\frac{w_2-\rho}{\Delta}%
\right) - T.
\end{equation*}
Therefore, the complete analytical solution for $k_2(t)$ is given by:
\begin{equation}
k_2(t)= \frac{\rho-w_2}{\beta} + \frac{\Delta}{\beta}\,\tanh\!\left[\frac{%
\Delta}{2}\left(t+\frac{2}{\Delta}\operatorname{arctanh}\!\left(\frac{w_2-\rho}{%
\Delta}\right)-T\right)\right], \, \forall t\in [0,T]. \label{eq:k2:append_del}
\end{equation}

Case 2: Suppose that $2\rho w_2 - \rho^2\geq 0$. Let $\Theta = \sqrt{2\rho w_2 -
\rho^2.}$ We can write the Riccati differential equation in \eqref{eq:riccati_st} as follows: $\dot{u}(t) = -(\beta/2)\left((u(t) - \rho/\beta)^2 +
(\Theta/\beta)^2\right).$
On rearranging, we obtain: 
\begin{equation*}
\frac{d u(t)}{(u(t) - \rho/\beta)^2 + (\Theta/\beta)^2} = -(\beta/2)
dt.
\end{equation*}
Next, we define $m(t):=\frac{u(t)-\rho/\beta}{\Theta/\beta}$. Then, we
obtain: $
d m(t)/(1+m(t)^2) = -(\Theta/2) dt.$
Taking integral on both sides, we obtain: $
\arctan \left(m(t) \right) = -(\Theta/2) (t+\zeta_\Theta).$
On substituting $m(t)=\frac{u(t)-\rho/\beta}{\Theta/\beta}$ 
and $u(t)=k_2(t)+\frac{w_2}{\beta} $, we obtain: 
\begin{equation*}
k_2(t) = \frac{\rho-w_2 + \Theta \tan\left(-\frac{\Theta}{2} (t +
\zeta_\Theta)\right)}{\beta}.
\end{equation*}
Using $k_2(T)=0$, we obtain: 
\begin{equation*}
\zeta_\Theta= -T-\frac{2}{\Theta}\operatorname{arctan}\left(\frac{w_2-\rho}{\Theta}%
\right).
\end{equation*}
 Therefore, the complete analytical solution for $k_2(t)$ is given by:
\begin{equation}
k_2(t)= \frac{\rho-w_2}{\beta} + \frac{\Theta}{\beta}\,\tan\!\left[-\frac{%
\Theta}{2}\left(t-\frac{2}{\Theta}\operatorname{arctan}\!\left(\frac{w_2-\rho}{%
\Theta}\right)-T\right)\right], \, \forall t\in [0,T]. \label{eq:k2:append_theta}
\end{equation}

\subsection{Proof for the analytical expression of $k_1(t)$}
\label{appen_k1} Recall that in the subsidy-active phase, the function $%
k_1(t)$ for $t\in (\tau_i, \tau_{i+1})$ satisfies 
\begin{equation}
\rho\,k_1(t)-\dot{k}_1(t)=\frac{\beta}{2}\Biggl(\frac{w_1}{\beta}+ k_1(t)+
s(\tau_i)+\eta_i\Biggr) \Biggl(\frac{w_2}{\beta}+ k_2(t)\Biggr),
\label{eq:k1_original}
\end{equation}
where $s(\tau_i)$ is subsidy level in the interval $(\tau_{i-1}, \tau_i]$ and $\eta_i$ is the subsidy adjustment at time $\tau_i$. Define $v_1(t)=w_1/\beta+ k_1(t)+ s(\tau_i)+\eta_i.$
Then, \eqref{eq:k1_original} can be rewritten as 
\begin{equation}
\rho\Bigl(v_1(t)- s(\tau_i)-\eta_i-w_1/\beta\Bigr)-\dot{v}_1(t)=\frac{\beta}{2%
}\,v_1(t) \left(\frac{w_2}{\beta}+ k_2(t)\right).\label{appen:v1}
\end{equation}
Case 1: Suppose that $\rho^2-2\rho w_2\geq 0$. Using the solution for $k_2(t)$
derived earlier in \eqref{eq:k2:append_del}, we obtain $\rho\Bigl(v_1(t)- s(\tau_i)-\eta_i-w_1/\beta\Bigr)-\dot{v}_1(t)=(1/2)%
\,v_1(t)[\rho+\Delta \tanh\!\left((\Delta/2)(t+\zeta_\Delta)%
\right)].$
Rearrange the terms to obtain %
\begin{equation}
\dot{v}_{1}(t)-\left[ \frac{\rho }{2}-\frac{\Delta }{2}\,\tanh \!\left( 
\frac{\Delta }{2}(t+\zeta _{\Delta })\right) \right] v_{1}(t)=-\rho \left(
s(\tau_i)+\eta_i+\frac{w_{1}}{\beta }\right) .  \label{eq:v1_ode}
\end{equation}%
The integrating factor is 
\begin{equation*}
I_{\Delta}(t)=\exp \!\left( -\int \left[ \frac{\rho }{2}-\frac{\Delta }{2}%
\,\tanh \!\left( \frac{\Delta }{2}(u+\zeta _{\Delta })\right) \right]
du\right) =\exp \!\left( -\frac{\rho }{2}\,t\right) \cosh \!\left( \frac{%
\Delta }{2}(t+\zeta _{\Delta })\right).
\end{equation*}%
Multiplying \eqref{eq:v1_ode} by $I_{\Delta }(t)$ and integrating from $t$
to $\tau _{i+1}$, we obtain 
\begin{equation*}
v_{1}(t)=\frac{1}{I_{\Delta }(t)}\left[ I_{\Delta }(\tau_{i+1}) v_{1}(\tau _{i+1})+\rho \left(
s(\tau_i)+\eta_i+\frac{w_{1}}{\beta }\right) \int_{t}^{\tau _{i+1}}I_{\Delta
}(u)\,du\right] .
\end{equation*}
Finally, returning to $k_{1}(t)$ via $
v_{1}(t)=k_{1}(t)+s(\tau_i)+\eta_i+w_{1}/\beta$,
we obtain the analytical expression for $k_{1}(t)$ for $t\in (\tau_i, \tau_{i+1}]$: 

\begin{align*}
k_{1}(t)&=\frac{1}{I_{\Delta }(t)}\left[ I_{\Delta }(\tau_{i+1}) \left(k_{1}(\tau _{i+1})+s(\tau
_{i})+\eta_i+w_{1}/\beta \right)+\rho \left( s(\tau _{i})+\eta_i+w_{1}/\beta \right) \int_{t}^{\tau _{i+1}}I_{\Delta }(u)\,du\right]\notag\\&\;\; -s(\tau
_{i})-\eta_i-w_{1}/\beta. 
\end{align*}%
Case 2: Suppose that $2\rho w_{2}-\rho ^{2}\geq 0$. Assume that, following our
model, we have introduced the shifted variable $v_{1}(t)=k_{1}(t)+s(\tau_i)+\eta_i+w_1/\beta$,
so that the dynamics of $v_{1}(t)$ in \eqref{appen:v1} after substituting \eqref{eq:k2:append_theta} satisfy the  differential equation: 
\begin{equation}
\dot{v}_{1}(t)-\left[ \frac{\rho }{2}-\frac{\Theta }{2}\,\tan \!\left( -%
\frac{\Theta }{2}(t+\zeta _{\Theta })\right) \right] v_{1}(t)=-\rho \left(
s(\tau_i)+\eta_i+w_{1}/\beta \right). \label{v1_del}
\end{equation}%
The integrating factor $I_\Theta(t)$ is given by 
\begin{equation*}
I_{\Theta }(t)=\exp \!\left( -\int \left[ \frac{\rho }{2}-\frac{\Theta }{2}%
\,\tan \!\left( -\frac{\Theta }{2}(u+\zeta _{\Theta })\right) \right]
du\right) =\exp \!\left( -\frac{\rho }{2}t\right)  \left| \cos \!\left( \frac{%
\Theta }{2}(t+\zeta _{\Theta })\right)\right|.
\end{equation*}%
Multiplying \eqref{v1_del} by $I_\Theta(t)$ and integrating from $t$ to $\tau_{i+1}$
yields 
\begin{equation*}
v_{1}(t)=\frac{1}{I_{\Theta }(t)}\left[ I_{\Theta }(\tau_{i+1})v_{1}(\tau_{i+1})+\rho \left(
s(\tau_i)+\eta_i+w_{1}/\beta \right) \int_{t}^{\tau _{i+1}}I_{\Theta
}(u)\,du\right] .
\end{equation*}%
Finally, reverting to $k_{1}(t)$ via $
v_{1}(t)=k_{1}(t)+s(\tau_i)+\eta_i+w_{1}/\beta,$ we obtain the solution for $k_{1}(t)$ for $t\in (\tau_i, \tau_{i+1}]$: 
\begin{align*}
k_{1}(t)=&\frac{1}{I_{\Theta }(t)}\left[ I_{\Theta }(\tau_{i+1})\left(k_{1}(\tau _{i+1})+s(\tau
_{i})+\eta_i+w_{1}/\beta\right)+\rho \left( s(\tau _{i})+\eta_i+\frac{w_{1}}{%
\beta }\right) \int_{t}^{\tau _{i+1}}I_{\Theta }(u)\,du\right] \notag\\&-s(\tau
_{i})-\eta_i-w_{1}/\beta.
\end{align*}
\setlength{\bibsep}{0pt}
\bibliography{biblio}
\bibliographystyle{apalike}
\end{document}